\documentclass[preprint,3p,times,review]{elsarticle}

\usepackage{amssymb}
\usepackage{amsmath}
\usepackage{amsthm}

\usepackage{lineno}

\DeclareFontFamily{U}{matha}{\hyphenchar\font45}
\DeclareFontShape{U}{matha}{m}{n}{
      <5> <6> <7> <8> <9> <10> gen * matha
      <10.95> matha10 <12> <14.4> <17.28> <20.74> <24.88> matha12
      }{}
\DeclareSymbolFont{matha}{U}{matha}{m}{n}
\DeclareFontSubstitution{U}{matha}{m}{n}

\DeclareFontFamily{U}{mathx}{\hyphenchar\font45}
\DeclareFontShape{U}{mathx}{m}{n}{
      <5> <6> <7> <8> <9> <10>
      <10.95> <12> <14.4> <17.28> <20.74> <24.88>
      mathx10
      }{}
\DeclareSymbolFont{mathx}{U}{mathx}{m}{n}
\DeclareFontSubstitution{U}{mathx}{m}{n}

\DeclareMathDelimiter{\vvvert}{0}{matha}{"7E}{mathx}{"17}

\newcommand{\bs}[1]{\boldsymbol{#1}}
\newcommand{\bbA}{\mathbb{A}}
\newcommand{\norm}[1]{\left\vvvert #1 \right\vvvert}
\newcommand{\calT}{\mathcal{T}}
\newcommand{\td}[1]{\tilde{#1}}

\theoremstyle{remark}
\newtheorem{remark}{Remark}[section]

\usepackage[overload]{empheq}

\usepackage{subcaption}

\usepackage{multirow}

\usepackage[linesnumbered,algoruled]{algorithm2e} 

\usepackage{enumitem}
\usepackage{booktabs}

\usepackage{hyperref}

\begin{document}

\begin{frontmatter}



\title{Certification of PGD reduced-order models with separated spatial variables}


\author[Navier]{Jean Ruel\corref{cor1}} 
\ead{jean.ruel@enpc.fr}
\cortext[cor1]{Corresponding author.}

\author[LMPS]{Ludovic Chamoin} 
\ead{ludovic.chamoin@ens-paris-saclay.fr}

\author[Navier,INRIA]{Frédéric Legoll}
\ead{frederic.legoll@enpc.fr}

\author[Navier]{Arthur Lebée}
\ead{arthur.lebee@enpc.fr}

\address[Navier]{Navier, ENPC, Institut Polytechnique de Paris, Université Gustave Eiffel, CNRS, Marne-la-Vallée, France}

\address[INRIA]{MATHERIALS project-team, Inria, Paris, France}

\address[LMPS]{Université Paris-Saclay, CentraleSupélec, ENS Paris-Saclay, CNRS, LMPS, Gif-sur-Yvette, France}

\begin{abstract}
Model order reduction techniques have become an attractive approach for obtaining fast approximations of multidimensional problems. Besides computational efficiency, ensuring the reliability of the resulting approximations is of primary importance. This work focuses on the certification of PGD-based reduced-order models based on the separation of spatial variables, which are particularly well suited to plate and shell geometries. Considering diffusion problems defined in plate-like domains, we introduce a guaranteed global error estimate associated with the PGD approximation. To this end, the error bounds are derived from the Constitutive Relation Error (CRE) method. The main difficulty of this approach lies in the construction of equilibrated fluxes, for which a dedicated procedure is proposed. Based on the resulting estimator, an adaptive strategy is developed to control both the discretization error and the number of PGD modes. This certification procedure is further extended to the error control in quantities of interest. We provide several numerical examples illustrating the reliability and efficiency of our procedure.
\end{abstract}



\begin{keyword}
Model Order Reduction \sep Proper Generalized Decomposition \sep Verification \sep Error estimation \sep Adaptivity 


\end{keyword}

\end{frontmatter}



\section{Introduction and objectives}

Numerical simulation has become an essential component of engineering analysis and design, and various numerical methods are available depending on the problem to be solved. Standard numerical methods can be limited when it comes to simulating multidimensional models in real time, and model order reduction techniques have thus emerged as an effective alternative over the past decades. They exploit the fact that the full-order solution of complex numerical models can often be accurately approximated by the reduced-order solution of surrogate models, so that the dimensionality can be drastically reduced. While enabling complex models to be solved in real time is an important challenge, ensuring the accuracy of the approximations obtained is equally important and can be crucial depending on the application. Among the most widely used approaches are Proper Orthogonal Decomposition (POD) and Reduced Basis (RB) methods, for which extensive work has been devoted to error estimation and certification.

In this work, we focus on the Proper Generalized Decomposition (PGD) method \cite{nouy_generalized_2008,nouy_priori_2010,chinesta_short_2011,chinesta_proper_2014}. Initially introduced as radial loading approximation \cite{ladeveze_nonlinear_1999}, it can be seen as a POD extension in which no \textit{a priori} knowledge on the solution is required. The PGD method is based on a modal representation of the solution with separation of variables (also referred to as a low-rank tensor approximation), and relies on an iterative strategy in which a set of simple problems are solved. Over the past years, many works have been devoted to the offline construction of PGD reduced models with applications to a wide range of parametrized problems, each parameter (related to spatial variables, time, material behavior, geometry, boundary or initial conditions,\dots) being seen as an extra-coordinate. Indeed, classical numerical methods based on brute force (grid-based) discretization rapidly lead to huge computational costs and storage requirements, as the number of degrees of freedom grows exponentially with respect to the number of dimensions of the problem; this is the so-called curse of dimensionality. Furthermore, in order to assess the accuracy of the constructed PGD solutions, some \textit{a posteriori} error estimation tools have been developed \cite{ammar_error_2010,ladeveze_verification_2011,de_almeida_basis_2013,allier_proper_2015,chamoin_posteriori_2017,reis_error_2020,reis_error_2020-1,reis_error_2021} using extensions of classical verification procedures used in finite element analysis. In particular, a method based on duality under the Constitutive Relation Error (CRE) concept was investigated in \cite{chamoin_posteriori_2017,chamoin_robust_2012} for linear elliptic and parabolic problems; it provides guaranteed error bounds accounting for all sources of error, for both global (energy norm) error and error in quantities of interest, as well as specific indicators on various error sources: truncation of the PGD approximation and discretization error of the underlying numerical technique.

In the present work, we address the verification of PGD approximations based on separated spatial variables. Such decompositions are particularly attractive in many applications \cite{bognet_advanced_2012,bognet_separated_2014,giner_proper_2013,vidal_proper_2013}, since an $n$-dimensional problem can be reformulated as a sequence of lower-dimensional problems defined on one-dimensional (or, in some cases, two-dimensional) domains. This substantially reduces computational costs and enables the efficient simulation of complex systems. A first verification strategy for this class of PGD approximations was proposed in \cite{nadal_separated_2015}. Based on recovery techniques, it provided a fast indicator of the discretization error and was mainly designed for mesh adaptation. However, deriving guaranteed error bounds that account for all error sources is not straightforward. In particular, the construction of equilibrated fields is more involved than in other PGD settings. The main objective of the present work is therefore to overcome this difficulty and to develop a fully guaranteed verification procedure for PGD approximations based on the separation of spatial variables. For the sake of completeness, we note that this question also arises for model reduction procedures based on the HiMod approach \cite{perotto_coupled_2014} which uses similar reduction concepts with separated spatial variables.

We focus on diffusion problems defined in plate-like domains. Within the CRE framework, we derive guaranteed \textit{a posteriori} error bounds for PGD approximations involving the separation of spatial variables. One of the main contribution of the article is the construction of equilibrated fluxes compatible with separated representations, which makes the application of CRE-based verification possible. The resulting estimator is then used to drive an adaptive strategy to control both the discretization error and the number of PGD modes. The proposed verification strategy is further extended to goal-oriented error estimation for quantities of interest. 

The article is organized as follows. After describing the reference problem, Section~\ref{sec:ref_problem} presents the basics of PGD model order reduction involving separation of spatial variables. Section~\ref{sec:err_est} is devoted to reviewing \textit{a posteriori} error estimation using CRE, as well as to defining and assessing a guaranteed estimator of the total PGD error. Based on this error estimator, an adaptive PGD strategy is proposed in Section~\ref{sec:adaptive_pgd} and illustrated numerically. Section~\ref{sec:goal_oriented} shows how to extend the previous tools to the goal-oriented error framework. Finally, conclusions and perspectives are outlined in Section~\ref{sec:conclusion}. 

\section{Reference problem and PGD model order reduction} \label{sec:ref_problem}

\subsection{Reference model and notations}

We consider the following $d$-dimensional diffusion problem, with in practice $d\in\{2,3\}$:
\begin{equation} \label{eq:ref_problem}
\begin{aligned}
    -\nabla\cdot (\bbA\bs{\nabla} u) &= f && \quad \textrm{in }\ \Omega = \omega\times I,\\
    \quad u &= 0 && \quad \textrm{on }\ \Gamma^D = \left(\Gamma^D_\omega\times I\right)\cup\left(\omega\times\Gamma^D_I\right),\\
    (\bbA\bs{\nabla} u)\cdot\bs{n} &= g && \quad \textrm{on } \ \Gamma^N = \left(\Gamma^N_\omega\times I\right)\cup\left(\omega\times\Gamma^N_I\right),
\end{aligned}
\end{equation}
where $\omega$ is an open bounded subset of $\mathbb{R}^{d-1}$ with boundary $\partial\omega$ and $\displaystyle I=\left(-\frac{t}{2},\frac{t}{2}\right)$ is an interval of $\mathbb{R}$ with boundary $\displaystyle\partial I=\left\{-\frac{t}{2},\frac{t}{2}\right\}$. The problem is defined on a plate-like domain. Furthermore, $\Gamma^D_\omega$ and $\Gamma^N_\omega$ are parts of $\partial\omega$ such that $\overline{\Gamma^D_\omega\cup\Gamma^N_\omega}=\partial\omega$ and $\Gamma^D_\omega\cap\Gamma^N_\omega=\emptyset$ (and likewise for $\Gamma^D_I$ and $\Gamma^N_I$) with $|\Gamma^D_\omega|\ne 0$ or $|\Gamma^D_I|\ne 0$. Boundary conditions are thus supposed to be compatible with separation of space variables. Observing that $\partial\Omega=(\partial\omega\times I)\cup(\omega\times\partial I)$, we define $\Gamma^D = \left(\Gamma^D_\omega\times I\right)\cup\left(\omega\times\Gamma^D_I\right)$ and $\Gamma^N = \left(\Gamma^N_\omega\times I\right)\cup\left(\omega\times\Gamma^N_I\right)$ as parts of $\partial\Omega$ such that $\overline{\Gamma^D\cup\Gamma^N}=\partial\Omega$, $\Gamma^D\cap\Gamma^N=\emptyset$ and $|\Gamma^D|\ne0$. We assume that $(f,g)\in L^2(\Omega)\times L^2(\Gamma^N)$ and that $\mathbb{A}\in\left[L^\infty(\Omega)\right]^{d\times d}$ is a symmetric positive definite matrix bounded from below. The flux associated with $u$ is denoted $\bs{q}=\bbA\bs{\nabla}u$ and $\bs{n}$ is the outer normal to $\partial \Omega$.

\smallskip 

Considering the functional space $V=\{v\in H^1(\Omega), \ v=0 \ \mathrm{on}\ \Gamma^D\}$, the weak formulation of \eqref{eq:ref_problem} consists in finding $u\in V$ such that,
\begin{equation} \label{eq:ref_weak_form}
    \forall\, v\in V, \quad B_1(u,v) = L\,(v),
\end{equation}
where
\begin{equation*}
    B_1(u,v) = \int_\Omega \bbA\bs{\nabla}u\cdot \bs{\nabla}v, \quad L\,(v) = \int_\Omega f\,v + \int_{\Gamma^N}g\,v.
\end{equation*}
This problem is equivalent to the following minimization problem:
\begin{equation*}
    u = \underset{v\in V}{\textrm{argmin }} J_1(v),
\end{equation*}
where $J_1$ is the potential energy functional defined by
\begin{equation*}
    \forall\, v\in V, \quad J_1(v) = \frac{1}{2}B_1(v,v) - L\,(v).
\end{equation*}
In what follows, we use the notation $\norm{\cdot}$ for the energy norm $\displaystyle \norm{v}=\sqrt{B_1(v,v)}$ on $V$ and, for any vector-valued field $\bs{p}\in \left(L^2(\Omega)\right)^d$, 
\begin{equation*}
    \norm{\bs{p}}_q = \sqrt{\int_\Omega \bbA^{-1}\bs{p}\cdot\bs{p}}.
\end{equation*}

\smallskip

The solution to \eqref{eq:ref_weak_form} is classically approximated using a finite element (FE) method. Introducing a FE space $V^h\subset V$, we recall that the FE approximation $u^h\in V^h$ of $u\in V$ is such that, 
\begin{equation*} 
    \forall\, v\in V^h, \quad B_1(u^h,v) = L\,(v).
\end{equation*}
In the next section, we introduce an alternative and computationally less expensive approximate solution of $u$ using PGD model order reduction. 

\subsection{PGD model order reduction}

We are interested here in approximating the solution to \eqref{eq:ref_weak_form} by a PGD strategy whose principle is to \textit{a priori} construct an approximation $u_m$ of $u$ as a separated variables representation defined in tensor product spaces, i.e. a finite sum of functions with separated variables. In the present context, involving separation of space variables, $u_m$ is of the form
\begin{equation} \label{eq:modal_dec}
    u_m(\bs{x},z) = \sum_{k=1}^{m} r_k(\bs{x})\, s_k(z), \quad (\bs{x},z)\in\omega\times I,
\end{equation}
where $m\in\mathbb{N}^*$ is the rank (or order) of the PGD approximation and, for any $1\leq k\leq m$, $(r_k,s_k)\in V_\omega\times V_I$  with $V_\omega=\{r\in H^1(\omega), \ r=0 \ \mathrm{on}\ \Gamma^D_\omega\}$ and $V_I=\{s\in H^1(I), \ s=0 \ \mathrm{on}\ \Gamma^D_I\}$. 

This decomposition is classically obtained using a greedy algorithm in which each term appearing in the sum \eqref{eq:modal_dec} is iteratively computed. At iteration $m$ of the algorithm, we look for a pair $(r_m,s_m)\in V_\omega\times V_I$ solution to the minimization problem
\begin{equation*}
    (r_m,s_m) \in \underset{(r,s)\,\in V_\omega\times V_I}{\textrm{argmin }} J_1(u_{m-1}+r\otimes s)
\end{equation*}
denoting $r \otimes s$ the tensor product $r \otimes s\,(\bs{x},z)=r\,(\bs{x})\,s\,(z)$ and where $\displaystyle u_{m-1}=\sum_{k=1}^{m-1}r_k \otimes s_k$ is the sum of terms computed at previous iterations. Assuming that a minimizer $(r_m,s_m)\in V_\omega\times V_I$ exists, it satisfies the following Euler-Lagrange equation:
\begin{equation} \label{eq:euler_lagrange}
    \forall\, (r,s)\in V_\omega\times V_I, \quad B_1(r_m\otimes s_m, r_m\otimes s + r\otimes s_m) = R_{m-1}(r_m\otimes s + r\otimes s_m),
\end{equation}
where 
\begin{equation*}
    R_{m-1}(v) = L\,(v) - B_1(u_{m-1},v)
\end{equation*}
is the residual at order $(m-1)$. Equation \eqref{eq:euler_lagrange} can be written equivalently as a system of coupled equations:
\begin{equation} \label{eq:pgd_problem}
\empheqlbrace\,
\begin{aligned}
    \forall\, r\in V_\omega, \quad B_1(r_m\otimes s_m, r\otimes s_m) &= R_{m-1}(r \otimes s_m),  \\ 
    \forall\, s\in V_I, \quad B_1(r_m\otimes s_m, r_m\otimes s) &= R_{m-1}(r_m\otimes s). 
\end{aligned}   
\end{equation}
From a numerical point of view, we introduce the discrete counterpart of \eqref{eq:pgd_problem}. Let $\calT_\omega$ (resp. $\calT_I$) be a partition of $\omega$ (resp. $I$). We denote $h_\omega$ and $h_I$ the maximum size of the elements of $\calT_\omega$ and $\calT_I$ and we introduce FE spaces $V_\omega^h\subset V_\omega$ and $V_I^h\subset V_I$ associated with $\calT_\omega$ and $\calT_I$. At iteration $m$, the PGD problem thus consists in finding $(r_m^h,s_m^h)\in V_\omega^h\times V_I^h$ such that,
\begin{subequations} \label{eq:pgd_problem_dis}
\begin{align}[left = \empheqlbrace\,]
    \forall\, r^h\in V_\omega^h, \quad B_1(r_m^h\otimes s_m^h, r^h\otimes s_m^h) &= R_{m-1}^h(r^h\otimes s_m^h), \label{eq:pgd_problem_1} \\ 
    \forall\, s^h\in V_I^h, \quad B_1(r_m^h\otimes s_m^h, r_m^h\otimes s^h) &= R_{m-1}^h(r_m^h\otimes s^h), \label{eq:pgd_problem_2} \
\end{align}
\end{subequations}
where
\begin{equation*}
    R^h_{m-1}(v) = L(v) - B_1(u^h_{m-1},v).
\end{equation*}
In practice, the system of equations \eqref{eq:pgd_problem_dis} is solved using a fixed-point algorithm. Initial \textit{ad hoc} functions $(r_m^{h,0},s_m^{h,0})$ are chosen. Then, at each step $i\geq1$, the algorithm computes $(r_{m}^{h,i},s_m^{h,i})$ such that
\begin{itemize}
    \item $r_m^{h,i}$ satisfies equation \eqref{eq:pgd_problem_1} for $s_m^h$ set to $s_m^{h,i-1}$;
    \item $s_m^{h,i}$ satisfies equation \eqref{eq:pgd_problem_2} for $r_m^h$ set to $r_m^{h,i}$.
\end{itemize}
The iterations stop when 
\begin{equation*}
    \frac{\norm{r_m^{h,i}\otimes s_m^{h,i} - r_m^{h,i-1}\otimes s_m^{h,i-1}}}{\norm{r_m^{h,i-1}\otimes s_m^{h,i-1}}} \leq \varepsilon_{FP}
\end{equation*}
where $\varepsilon_{FP}$ is a predefined tolerance threshold (where the subscript FP stands for fixed-point). A maximum number of iterations is also imposed to ensure that the fixed-point algorithm terminates. 

\medskip
At the end, the computed solution is
\begin{equation*}
    u_m^{\bs{h}} = \sum_{k=1}^m r_k^h \otimes s_k^h
\end{equation*}
with $\bs{h}=(h_\omega,h_I)$. The overall error of the PGD approach is defined by $e_{PGD}=u-u^{\bs{h}}_m$. The purpose of Section~\ref{sec:err_est} is to design an estimator in order to assess this error. At this stage, we can already mention that the error in the PGD approach stems from two main causes: 
\begin{itemize}
    \item the discretization, due to the introduction of FE spaces $V_\omega^h$ and $V_I^h$ to compute each PGD mode;
    \item the reduction, due to the truncation at mode $m$ of the modal decomposition, and to the stopping criterion used in the fixed-point procedure related to \eqref{eq:pgd_problem_dis}.
\end{itemize}

\begin{remark}
In this article, we do not consider other error sources than those arising from discretization and reduction. Other potential error sources (including quadrature error or algebraic error due to the use of iterative solvers for linear systems) are assumed to be controlled and negligible in comparison. \qed
\end{remark}

\section{Global error estimation for the PGD approximation} \label{sec:err_est}

Our aim is to provide tools for estimating the error $e_{PGD}$ between the exact solution $u$ and the approximate solution $u^{\bs{h}}_m$ computed by PGD. For this purpose, we define and assess in this section a guaranteed \textit{a posteriori} estimator of the global error $e_{PGD}$ measured in the energy norm. 

\subsection{\textit{A posteriori} error estimation using CRE}

In order to construct a guaranteed error estimator, we resort to the Constitutive Relation Error (CRE) concept which leads to define a dual formulation of Problem~\eqref{eq:ref_problem} and the space of equilibrated fluxes
\begin{equation*}
    W = \left\{\bs{p}\in H(\textrm{div},\Omega),\ \bs{\nabla}\cdot\bs{p}+f=0\textrm{ in }\Omega,\ \bs{p}\cdot\bs{n}=g\textrm{ on }\Gamma^N\right\},
\end{equation*}
where $\displaystyle H(\textrm{div},\Omega)=\left\{\bs{p}\in\left[L^2(\Omega)\right]^d,\ \bs{\nabla}\cdot\bs{p}\in L^2(\Omega)\right\}$. A flux field that belongs to $W$ is said to be statically admissible (SA). In the same way, any $\hat{u}\in V$ is said to be kinematically admissible (KA), which is the case for $u^{\bs{h}}_m$. Then, for the admissible pair $(u_m^{\bs{h}},\hat{\bs{q}})\in V\times W$, the CRE functional $E_{\text{CRE}}$ is defined by
\begin{equation*}
E_{\text{CRE}}\left(u_m^{\bs{h}},\hat{\bs{q}}\right) = \norm{\hat{\bs{q}}-\bbA\bs{\nabla}u_m^{\bs{h}}}_q. 
\end{equation*}
Introducing the complementary energy functional $J_2$ defined by
\begin{equation} \label{eq:comp_energy}
    \forall\, \bs{p}\in \left(L^2(\Omega)\right)^d, \quad J_2(\bs{p}) = \frac{1}{2}\int_\Omega \bbA^{-1}\bs{p}\cdot\bs{p},
\end{equation}
we also have
\begin{equation*}
    E_{\text{CRE}}^2\left(u_m^{\bs{h}},\hat{\bs{q}}\right) = 2\left(J_1(u_m^{\bs{h}})+J_2(\hat{\bs{q}})\right),
\end{equation*}
which follows from the fact that
\begin{equation} \label{eq:weak_equilibrium}
    \forall\, v\in V,\quad \int_\Omega \hat{\bs{q}}\cdot \bs{\nabla}v = \int_\Omega f\,v +\int_{\Gamma^N} g\, v.
\end{equation}
Recalling that the exact flux field $\bs{q}=\bbA\bs{\nabla}u$ is the unique minimizer of the complementary problem
\begin{equation} \label{eq:complementary_problem}
    \bs{q} = \underset{\bs{p}\in W}{\textrm{argmin }} J_2(\bs{p}),
\end{equation}
it is not difficult to show that
\begin{equation} \label{eq:CRE_upper_bound}
    \forall\,\hat{\bs{q}}\in W,\quad \norm{u-u_m^{\bs{h}}}^2 = E_{\text{CRE}}^2\left(u_m^{\bs{h}},\bs{q}\right) \leq E_{\text{CRE}}^2\left(u_m^{\bs{h}},\hat{\bs{q}}\right)=2\left(J_1(u_m^{\bs{h}})+J_2(\hat{\bs{q}})\right).
\end{equation}
Another way of presenting \eqref{eq:CRE_upper_bound} is the so-called Prager-Synge equality, which is
\begin{equation} \label{eq:prager_synge}
    \forall\,\hat{\bs{q}}\in W,\quad E_{\text{CRE}}^2\left(u_m^{\bs{h}},\hat{\bs{q}}\right) = \norm{u-u_m^{\bs{h}}}^2 + \norm{\bs{q}-\hat{\bs{q}}}_q^2.
\end{equation}
From this equality follows the hypercircle property, which will be useful in Section~\ref{sec:goal_oriented}:
\begin{equation} \label{eq:hypercircle}
    \forall\,\hat{\bs{q}}\in W,\quad E_{\text{CRE}}\left(u_m^{\bs{h}},\hat{\bs{q}}\right) = 2 \norm{\bs{q} - \bs{q}^*}_q, \quad \text{where} \quad \bs{q}^* = \frac{1}{2}\left(\hat{\bs{q}} + \bbA \bs{\nabla}u_m^{\bs{h}}\right).
\end{equation}
In view of \eqref{eq:CRE_upper_bound} or \eqref{eq:prager_synge}, for any $\hat{\bs{q}}\in W$, $E_{\text{CRE}}\left(u_m^{\bs{h}},\hat{\bs{q}}\right)$ provides an upper bound on the PGD error measured in the energy norm. Subject to the construction of an appropriate $\hat{\bs{q}}\in W$, this quantity constitutes our error estimator. 

\subsection{Complementary PGD approach} \label{sec:complementary_PGD}

The challenge with the CRE concept lies in the construction of a relevant statically admissible flux $\hat{\bs{q}}\in W$, that is, a flux sufficiently close to $\bs{q}$ in order to obtain an accurate error estimator (see \eqref{eq:prager_synge}). In the FEM context, this flux field is generally obtained by post-processing the FE flux $\bs{q}^h=\bbA\bs{\nabla}u^h$ and solving local problems at the element scale (see e.g. \cite{chamoin_introductory_2022} for a recent review). The methods using this strategy are largely based on the fact that $\bs{q}^h$ satisfies the weak equilibrium \eqref{eq:weak_equilibrium} in the FE sense, i.e. for any $v\in V^h$. However, the PGD flux $\bs{q}^{\bs{h}}_m=\bbA\bs{\nabla}u^{\bs{h}}_m$ with separated spatial variables does not satisfy this weak equilibrium, which makes it impossible to directly use existing techniques to reconstruct a SA flux. 

In what follows, we proceed differently. The PGD method is used to approximate the complementary problem \eqref{eq:complementary_problem}. We therefore define a $n$-order PGD approximation $\hat{\bs{q}}_n$ of the flux as follows:
\begin{equation} \label{eq:pgd_flux_form_1}
    \hat{\bs{q}}_{n}(\bs{x},z) = \bs{q}_0(\bs{x},z) + \sum_{k=1}^{n}\bs{\Phi}_k(\bs{x})\circ \bs{\Psi}_k(z),
\end{equation}
where $\circ$ denotes the elementwise product. The term $\bs{q}_0$ is a particular flux that belongs to $W$, and, for any $1\leq k\leq n$, $\bs{\Phi}_k\circ \bs{\Psi}_k$ is statically admissible to zero, that is
\begin{equation} \label{eq:div0_const}
    \nabla\cdot(\bs{\Phi}_k\circ \bs{\Psi}_k) = 0\textrm{ in } \Omega, \quad (\bs{\Phi}_k\circ \bs{\Psi}_k)\cdot\bs{n}=0 \textrm{ on } \Gamma^N.
\end{equation}

\begin{remark}
The question arises as to the choice of $n$ for the PGD approximation of the flux. This choice has an influence on the quality of the resulting estimator. By default, we have chosen $n=m$ where $m$ is the number of modes in \eqref{eq:modal_dec}. This question is discussed from a numerical viewpoint in Section~\ref{sec:err_est_test}. \qed
\end{remark}

In order to satisfy \eqref{eq:div0_const}, the form of each term $\bs{\Phi}_k\circ \bs{\Psi}_k$ is imposed \textit{a priori} in a separated form:
\begin{equation} \label{eq:pgd_flux_form_2}
    \bs{\Phi}_k(\bs{x})\circ \bs{\Psi}_k(z) = 
    \begin{pmatrix}
        \bs{\varphi}_k(\bs{x})\, \psi_k'(z)\\
        -\left(\nabla\cdot \bs{\varphi}_k\right)(\bs{x})\, \psi_k(z)
    \end{pmatrix}
\end{equation}
with $\bs{\varphi}_k\cdot\bs{n}=0$ on $\Gamma_\omega^N$ and $\psi_k=0$ on $\Gamma_I^N$. In \eqref{eq:pgd_flux_form_2}, the prime denotes the derivative and $\bs{\varphi}_k$ is a vector of dimension $d-1$, while $\psi_k$ is scalar-valued (beware in \eqref{eq:pgd_flux_form_2} of the difference between $\bs{\Psi}_k$ and $\psi_k$).

With regard to $\bs{q}_0$, we can write $\bs{q}_0 = \bs{q}_0^f + \bs{q}_0^g$ where $\bs{q}_0^f$ and $\bs{q}_0^g$ are such that 
\begin{equation*}
\left\{
\begin{aligned}
    \nabla\cdot \bs{q}_0^f + f &= 0 && \quad \textrm{in }\ \Omega\\
    \bs{q}_0^f\cdot\bs{n} &= 0 && \quad \textrm{on } \ \Gamma^N
\end{aligned}
\right. \quad \textrm{and } 
\left\{
\begin{aligned}
    \nabla\cdot \bs{q}_0^g &= 0 && \quad \textrm{in }\ \Omega\\
    \bs{q}_0^g\cdot\bs{n} &= g && \quad \textrm{on } \ \Gamma^N
\end{aligned}
\right. .
\end{equation*}
In practice, $\bs{q}_0^f$ and $\bs{q}_0^g$ can be defined analytically in a straightforward manner for simple loading. In the representative case where $\displaystyle \Gamma^N = \omega\times \left\{\frac{t}{2}\right\}$ and assuming that $f$ can be written as a sum of functions with separate variables, that is 
\begin{equation*}
    f(\bs{x},z) = \sum_{k=1}^{m_f} f_k^\omega(\bs{x})\,f_k^I(z),
\end{equation*}
$\bs{q}_0^f$ can be defined as
\begin{equation} \label{eq:q0f}
    \bs{q}_0^f(\bs{x},z) = -\sum_{k=1}^{m_f}
    \begin{pmatrix}
        \bs{0} \\
        f_k^\omega(\bs{x})\,F_k^I(z) 
    \end{pmatrix}
\end{equation}
where $\left(F_k^I\right)'(z) = f_k^I(z)$ and $\displaystyle F_k^I\left(\frac{t}{2}\right) = 0$.  We then define $\bs{q}_0^g$ by
\begin{equation} \label{eq:q0g}
    \bs{q}_0^g(\bs{x},z) = 
    \begin{pmatrix} 
        \bs{0} \\
        \displaystyle g(\bs{x})
    \end{pmatrix}.
\end{equation}
Note that $\bs{q}_0^g=\bs{0}$ when the problem \eqref{eq:ref_problem} is submitted to pure Dirichlet boundary conditions.

Motivated by \eqref{eq:CRE_upper_bound}, we now minimize the complementary energy \eqref{eq:comp_energy} on the flux of the form \eqref{eq:pgd_flux_form_1}-\eqref{eq:pgd_flux_form_2} using the PGD standard procedure. Let us define $W_\omega=\{\bs{\varphi}\in H(\textrm{div},\omega), \ \bs{\varphi}\cdot\bs{n}=0 \ \mathrm{on}\ \Gamma^N_\omega\}$ and $W_I=\{\psi\in H^1(I), \ \psi=0 \ \mathrm{on}\ \Gamma^N_I\}$. At iteration $n$, we look for a pair $(\bs{\varphi}_n,\psi_n)\in W_\omega\times W_I$ solution to the minimization problem
\begin{equation*}
    (\bs{\varphi}_n,\psi_n) \in \underset{(\bs{\varphi},\psi)\,\in W_\omega\times W_I}{\textrm{argmin }} J_2(\hat{\bs{q}}_{n-1}+\bs{\Phi}\circ\bs{\Psi})
\end{equation*}
where $\displaystyle \hat{\bs{q}}_{n-1}=\bs{q}_0 + \sum_{k=1}^{n-1}\bs{\Phi}_k\circ\bs{\Psi}_k$ is the sum of terms computed at previous iterations. Assuming that a minimizer $(\bs{\varphi}_n,\psi_n)\in W_\omega\times W_I$ exists, it satisfies the following Euler-Lagrange equation:
\begin{equation} \label{eq:dual_euler_lagrange}
    \forall\, (\bs{\varphi},\psi)\in W_\omega\times W_I, \quad B_2(\bs{\Phi}_n\circ\bs{\Psi}_n, \bs{\Phi}_n\circ\bs{\Psi} + \bs{\Phi}\circ\bs{\Psi}_n) = R_{n-1}(\bs{\Phi}_n\circ\bs{\Psi} + \bs{\Phi}\circ\bs{\Psi}_n),
\end{equation}
where 
\begin{equation*}
    B_2(\bs{q},\bs{p}) = \int_\Omega \bbA^{-1} \bs{q}\cdot \bs{p} \quad \textrm{and} \quad R_{n-1}(\bs{p}) = - B_2(\hat{\bs{q}}_{n-1},\bs{p})
\end{equation*}
is the residual at order $(n-1)$. Equation \eqref{eq:dual_euler_lagrange} can be written equivalently as a system of coupled equations:
\begin{subequations} \label{eq:dual_pgd_problem}
\begin{align}[left = \empheqlbrace\,]
    \forall\, \bs{\varphi}\in W_\omega, \quad B_2(\bs{\Phi}_n\circ\bs{\Psi}_n, \bs{\Phi}\circ\bs{\Psi}_n) &= R_{n-1}(\bs{\Phi}\circ\bs{\Psi}_n),  \\ 
    \forall\, \psi\in W_I, \quad B_2(\bs{\Phi}_n\circ\bs{\Psi}_n, \bs{\Phi}_n\circ\bs{\Psi}) &= R_{n-1}(\bs{\Phi}_n\circ\bs{\Psi}). 
\end{align}
\end{subequations}
From a numerical point of view, a discrete counterpart of \eqref{eq:dual_pgd_problem} is used. We perform here a finite element discretization for each subproblem. Note that $C^0$ continuous elements are sufficient to ensure normal flux continuity. 

\begin{remark}
    With regard to the discretization of the complementary PGD approach \eqref{eq:dual_pgd_problem}, standard Lagrange finite elements are used in all the numerical tests reported here, regardless of the dimension of the problem. When $d=2$, this choice is natural since $W_\omega$ is in fact defined by $W_\omega=\{\varphi\in H^1(\omega), \ \varphi=0 \ \mathrm{on}\ \Gamma^N_\omega\}$. When $d\geq3$, however, this choice implies that $\bs{\varphi}_k^h\in \left[H^1(\omega)\right]^{d-1}$, which is sufficient for $\bs{\varphi}_k^h$ to belong to $H(\textrm{div},\omega)$ but not necessary. \qed 
\end{remark}

\begin{remark}
    In general cases where $\displaystyle \Gamma^N \ne \omega\times \left\{\frac{t}{2}\right\}$, determining $\bs{q}_0$ may require a little more work. Let us assume that the flux boundary conditions are written as 
    \begin{equation*}
        \bs{q}\cdot\bs{n}\,(\bs{x},z) = \left\{
    \begin{array}{cl}
        \displaystyle g_t(\bs{x}) & \textrm{on } \omega\times\left\{\frac{t}{2} \right\}\\
        \displaystyle g_b(\bs{x}) & \textrm{on } \omega\times\left\{-\frac{t}{2} \right\}\\
        g_l(z) & \textrm{on } \Gamma_\omega^N\times I,
    \end{array}
    \right.
    \end{equation*}
    where we underline that the function $g_l$ does not depend on $\bs{x}$. On the one hand, once we have computed $\bs{\varphi}^f_k$ such that $\nabla\cdot \bs{\varphi}^f_k + f^\omega_k=0$ in $\omega$ and $\bs{\varphi}^f_k\cdot\bs{n}=0$ on $\Gamma^N_\omega$, the flux $\bs{q}_0^f$ can be defined as 
    \begin{equation*}
        \bs{q}_0^f(\bs{x},z) = \sum_{k=1}^{m_f}
        \begin{pmatrix}
            \bs{\varphi}^f_k(\bs{x})f^I_k(z) \\
            0
        \end{pmatrix}.
    \end{equation*}
    On the other hand, once $\bs{\varphi}^g_t$, $\bs{\varphi}^g_b$ and $\bs{\varphi}^g_l$ have been computed such that $\displaystyle \nabla\cdot \bs{\varphi}^g_t + \frac{g_t}{t}=0$ in $\omega$ and $\bs{\varphi}^g_t\cdot\bs{n}=0$ on $\Gamma^N_\omega$, likewise for $\bs{\varphi}^g_b$, and $\nabla\cdot \bs{\varphi}^g_l=0$ in $\omega$ and $\bs{\varphi}^g_l\cdot\bs{n}=1$ on $\Gamma^N_\omega$, $\bs{q}_0^g$ can be defined as 
    \begin{equation*}
        \bs{q}_0^g(\bs{x},z) =
        \begin{pmatrix}
            \bs{\varphi}^g_t(\bs{x}) + \bs{\varphi}^g_b(\bs{x}) + \bs{\varphi}^g_l(\bs{x})g_l(z) \\
            \displaystyle \left(\frac{z}{t} + \frac{1}{2}\right)g_t(\bs{x}) + \left(\frac{z}{t} - \frac{1}{2}\right)g_b(\bs{x})
        \end{pmatrix}.
    \end{equation*}
    The functions $\bs{\varphi}^f_k$, $\bs{\varphi}^g_t$, $\bs{\varphi}^g_b$, $\bs{\varphi}^g_l$ can be determined by duality using a standard \textit{a posteriori} equilibrated finite element computation. It should be noted that these computations are only in dimension $d-1$, and thus expected to be inexpensive. \qed
\end{remark}

\subsection{Numerical assessment of the global error estimator} \label{sec:err_est_test}

Once $u_m^{\bs{h}}$ and $\hat{\bs{q}}_n^{\bs{h}}$ have been computed, we define our global PGD error estimator $\eta_{mn}^{\bs{h}}$ by
\begin{equation} \label{eq:error_estimate}
    \eta_{mn}^{\bs{h}} = E_{\text{CRE}}\left(u_m^{\bs{h}},\hat{\bs{q}}_n^{\bs{h}}\right),
\end{equation}
as well as the effectivity index $I_{\text{eff}}$ by
\begin{equation*}
    I_{\text{eff}} = \frac{\eta_{mn}^{\bs{h}}}{\norm{u - u^{\bs{h}}_m}}.
\end{equation*}
In this section, we assess the accuracy of this estimator through two 2D test cases where an analytical solution is known. In particular, we solve the Poisson problem $-\Delta u = f$ defined in the unit square domain $\Omega = (0,1)^2$. We also consider the same problem posed on a thin domain in Section~\ref{sec:thin_domains}.

\subsubsection{2D Poisson problem with a uniform source term}

We start by considering homogeneous Dirichlet boundary conditions on $\partial\Omega$  and a uniform source term $f=1$. Consequently, according to \eqref{eq:q0f}, we choose
\begin{equation*}
    \bs{q}_0(x,z) =
    \begin{pmatrix}
    0 \\
    -z
    \end{pmatrix}
\end{equation*}
as a particular flux that equilibrates the external loading. The analytical solution to this problem is given by
\begin{equation*}
    u(x,z) = \frac{x\,(1-x)}{2} - \frac{4}{\pi^3} \sum_{\mathrm{odd}\ n=1}^\infty \sin{(n\pi x)}\frac{\sinh{(n\pi z)} + \sinh{(n\pi\,(1-z))}}{n^3\sinh{(n\pi)}}.
\end{equation*}
We note in passing that the separate-variable form of the exact solution conforts the use of PGD in this context.

\begin{figure}[h!]
    \centering
    \begin{subfigure}{0.49\textwidth}
        \includegraphics[trim=3cm 10cm 4cm 9cm, clip=true, width=\textwidth]{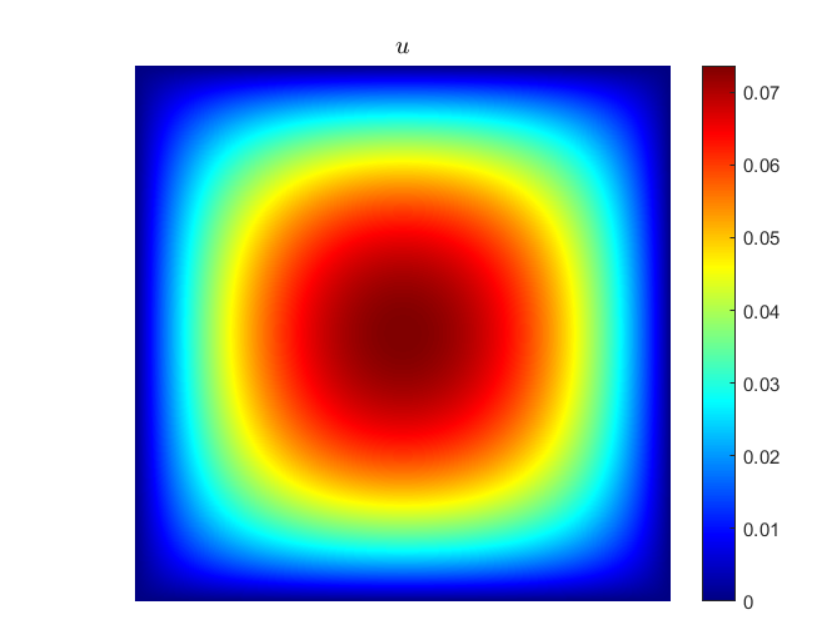}
    \end{subfigure}
    \begin{subfigure}{0.49\textwidth}
        \includegraphics[trim=3cm 10cm 4cm 9cm, clip=true, width=\textwidth]{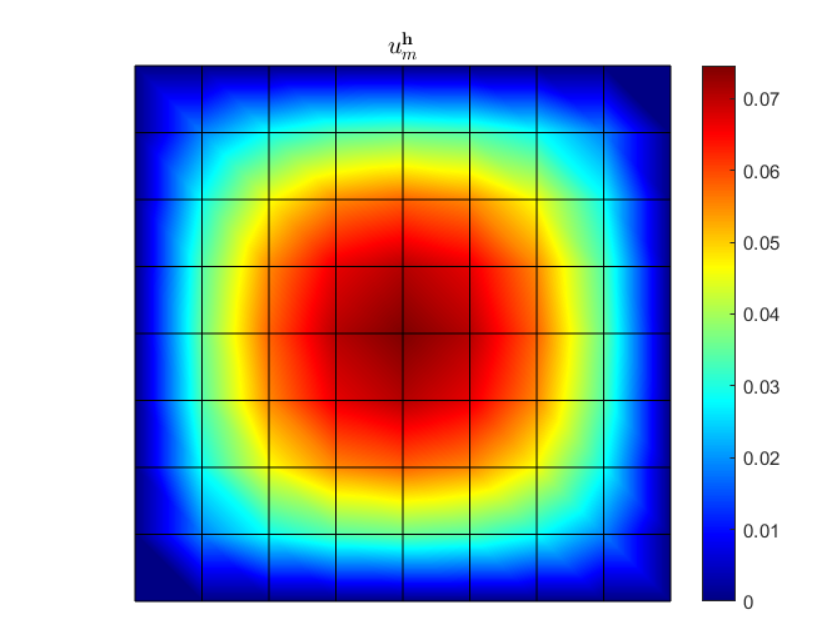}
    \end{subfigure}
    \caption{2D Poisson problem with a uniform source term: exact solution $u$ (left) and PGD solution with 3 modes (right).}
    \label{fig:sol_1}
\end{figure}

\begin{figure}[h!]
    \centering
    \begin{subfigure}{0.49\textwidth}
        \includegraphics[trim=3cm 10cm 4cm 9cm, clip=true, width=\textwidth]{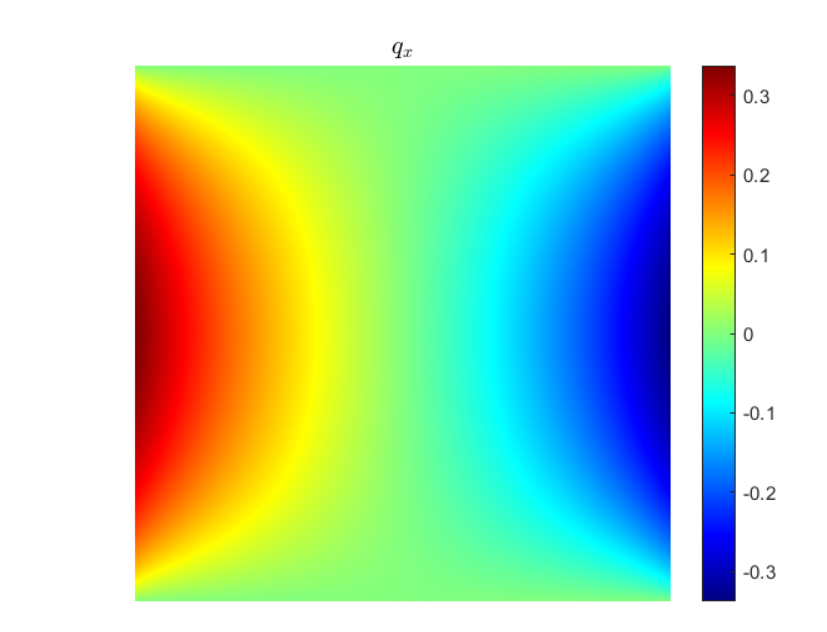}
    \end{subfigure}
    \begin{subfigure}{0.49\textwidth}
        \includegraphics[trim=3cm 10cm 4cm 9cm, clip=true, width=\textwidth]{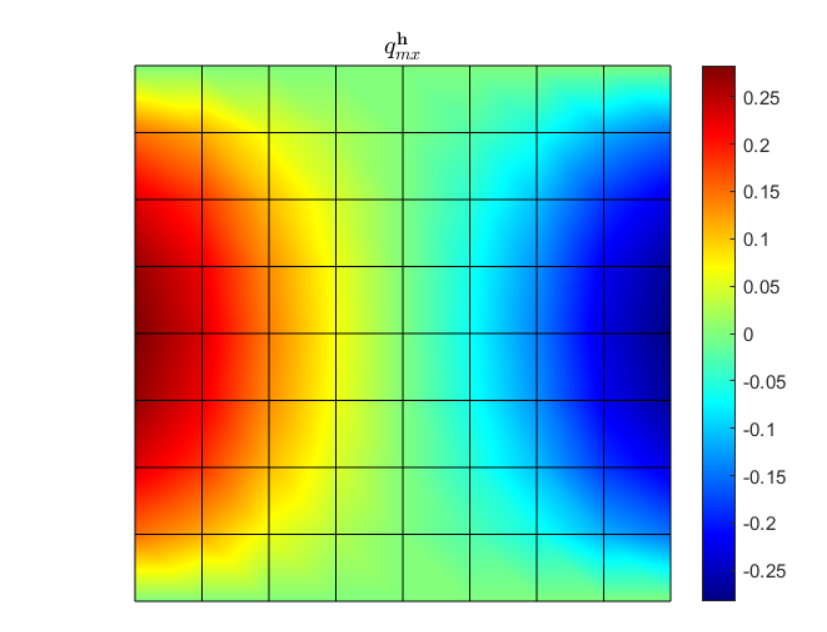}
    \end{subfigure}
    \vfill
    \begin{subfigure}{0.49\textwidth}
        \includegraphics[trim=3cm 10cm 4cm 9cm, clip=true, width=\textwidth]{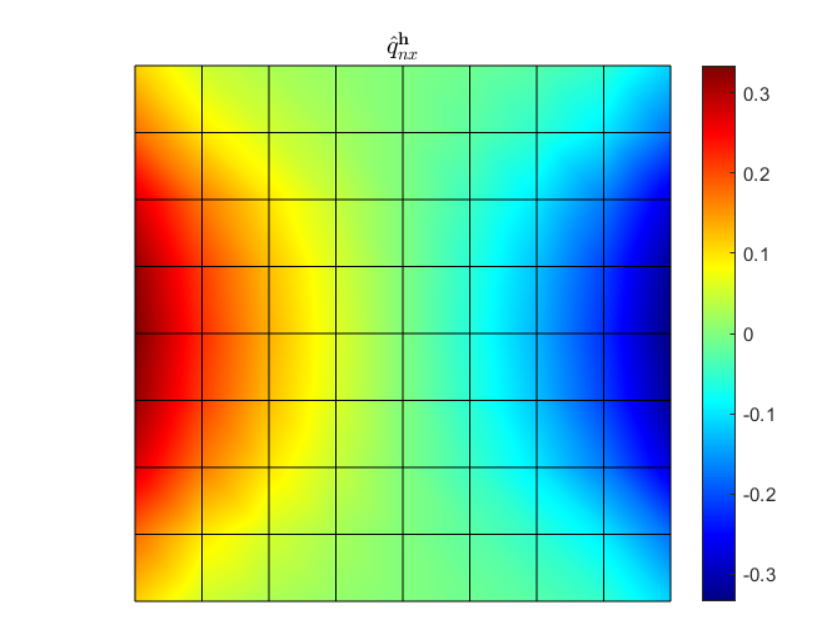}
    \end{subfigure}
    \caption{2D Poisson problem with a uniform source term: exact flux $\bs{q}$ (top left), PGD flux $\bs{q}^{\bs{h}}_m$ (top right) and SA PGD flux $\hat{\bs{q}}^{\bs{h}}_n$ with $n=m=3$ modes (bottom). We only show the first component of the flux; the figure for the second component is similar.}
    \label{fig:fluxes_x_1}
\end{figure}

An approximate PGD solution is computed using uniform meshes $\calT_\omega$ and $\calT_I$ composed of 8 piecewise linear elements for both PGD approaches. The PGD solution $u_m^{\bs{h}}$ is represented in Figure~\ref{fig:sol_1} with 3 modes and compared with the exact solution $u$. As suggested by the analytical form of the solution, only a few PGD modes are needed to recover the exact solution qualitatively and we observe a good agreement between $u$ and $u^{\bs{h}}_m$. The exact flux $\bs{q}$ and the PGD fluxes $\bs{q}^{\bs{h}}_m$ and $\hat{\bs{q}}^{\bs{h}}_n$ are shown in Figure~\ref{fig:fluxes_x_1}. On these visualizations, the reconstructed 2D mesh obtained from the mesh of $\omega$ and of $I$ is represented in black solid lines. Here again, both PGD fluxes are correctly approximated with only three modes. We simply note a deviation of the SA PGD flux $\hat{\bs{q}}^{\bs{h}}_n$ from the exact flux $\bs{q}$ at the corners of the domain. This flux does not satisfy any compatibility equation, since it does not derive from a gradient and is therefore not strictly zero at the top and bottom boundaries of the domain.
 
Regarding the accuracy of the error estimator \eqref{eq:error_estimate}, its effectivity index is shown in Table~\ref{tab:eff_index_1_1} with respect to $m$ for different choices of the rank $n$ of the complementary PGD approximation. In general, the effectivity indices are quite acceptable in this case and the error estimator is well guaranteed. Furthermore, it appears that it is preferable to take $n>m$, particularly for accurately estimating the PGD error when $m$ is very small. In Table~\ref{tab:eff_index_1_2}, we study the influence of discretization on the accuracy of the error estimator. We denote by $N_\omega$ (resp. $N_I$) the number of elements of $\calT_\omega$ (resp. $\calT_I$) and by $p$ the degree of elements used to compute $\hat{\bs{q}}^{\bs{h}}_n$. Thus, $\hat{\bs{q}}_n^{\bs{h}}$ is computed on the same meshes as $u_m^{\bs{h}}$ but with piecewise linear ($p=1$) or quadratic ($p=2$) elements, whereas only linear elements are considered for $u_m^{\bs{h}}$. We can see that the finer the meshes, the less accurate the error estimator is, although this trend is much less pronounced when $m$ increases. This can be explained by the fact that when the mesh is refined, the reduction error dominates and the small number of modes used to represented the SA PGD flux becomes apparent. On the other hand, computing $\hat{\bs{q}}_n^{\bs{h}}$ with elements of higher degree than for $u$ improves the accuracy of the error estimator, and does so significantly when $m$ is sufficiently large. In view of these results, we suggest computing the SA PGD flux with $n>m$ and with elements of higher degree $p$ than for the approximate PGD solution $u_m^{\bs{h}}$. The precise choice of $n$ and $p$ results from a balance between precision and computation time. For applications where it is expected that only a few PGD modes will be needed, $n=m+1$ is an acceptable choice (Table~\ref{tab:eff_index_1_1} shows that considering $n=m+2$ does not improve the results compared to the choice $n=m+1$).

\begin{table}[h!]
    \centering
    \begin{tabular}{cc|cc|cc|cc}
    \toprule
    $m$ & $\norm{u-u_m^{\bs{h}}}$ & $\eta_{mm}^{\bs{h}}$ & $I_{\text{eff}}$ & $\eta_{m\,(m+1)}^{\bs{h}}$ & $I_{\text{eff}}$ & $\eta_{m\,(m+2)}^{\bs{h}}$ & $I_{\text{eff}}$\\
    \midrule
    1 & $2.89\times 10^{-2}$ & $6.18\times 10^{-2}$ & 2.14 & $3.66\times 10^{-2}$ & 1.27 & $3.63\times 10^{-2}$ & 1.26 \\
    2 & $2.85\times 10^{-2}$ & $3.63\times 10^{-2}$ & 1.27 & $3.60\times 10^{-2}$ & 1.26 & $3.59\times 10^{-2}$ & 1.26 \\
    3 & $2.85\times 10^{-2}$ & $3.60\times 10^{-2}$ & 1.26 & $3.59\times 10^{-2}$ & 1.26 & $3.59\times 10^{-2}$ & 1.26 \\
    4 & $2.85\times 10^{-2}$ & $3.59\times 10^{-2}$ & 1.26 & $3.59\times 10^{-2}$ & 1.26 & $3.59\times 10^{-2}$ & 1.26 \\
    5 & $2.85\times 10^{-2}$ & $3.59\times 10^{-2}$ & 1.26 & $3.59\times 10^{-2}$ & 1.26 & $3.59\times 10^{-2}$ & 1.26 \\
    \bottomrule
    \end{tabular}
    \caption{2D Poisson problem with a uniform source term: effectivity index ($I_{\text{eff}}$) of the estimator $\eta_{mn}^{\bs{h}}$ with respect to the rank of the complementary PGD solution. In this example, $\norm{u}=1.88\times 10^{-1}$.}
    \label{tab:eff_index_1_1}
\end{table}

\begin{table}[h!]
    \centering
    \begin{tabular}{cc|ccc|ccc}
    \toprule
    $p$ & $N_\omega\times N_I$ & $\norm{u-u^{\bs{h}}_1}$ & $\eta_{11}^{\bs{h}}$ & $I_{\text{eff}}$ & $\norm{u-u^{\bs{h}}_4}$ & $\eta_{44}^{\bs{h}}$ & $I_{\text{eff}}$ \\
    \midrule
    \multirow{3}{*}{1}
    & $8\times 8$ & $2.89\times 10^{-2}$ & $6.18\times 10^{-2}$ & 2.14 & $2.85\times 10^{-2}$ & $3.59\times 10^{-2}$ & 1.26 \\
    & $16\times 16$ & $1.54\times 10^{-2}$ & $5.48\times 10^{-2}$ & 3.55 & $1.43\times 10^{-2}$ & $1.81\times 10^{-2}$ & 1.27 \\
    & $32\times 32$ & $9.40\times 10^{-3}$ & $5.29\times 10^{-2}$ & 5.63 & $7.15\times 10^{-3}$ & $9.14\times 10^{-3}$ & 1.28 \\
    \midrule
    \multirow{3}{*}{2}
    & $8\times 8$ & $2.89\times 10^{-2}$ & $5.94\times 10^{-2}$ & 2.06 & $2.85\times 10^{-2}$ & $2.86\times 10^{-2}$ & 1.00 \\
    & $16\times 16$ & $1.54\times 10^{-2}$ & $5.41\times 10^{-2}$ & 3.51 & $1.43\times 10^{-2}$ & $1.44\times 10^{-2}$ & 1.01 \\
    & $32\times 32$ & $9.40\times 10^{-3}$ & $5.27\times 10^{-2}$ & 5.61 & $7.15\times 10^{-3}$ & $7.30\times 10^{-3}$ & 1.02 \\
    \bottomrule
    \end{tabular}
    \caption{2D Poisson problem with a uniform source term: influence of the complementary PGD discretization (through the polynomial degree $p$) on the effectivity index with respect to the number of PGD modes. We recall that $\norm{u}=1.88\times 10^{-1}$.}
    \label{tab:eff_index_1_2}
\end{table}

\subsubsection{2D Poisson problem with a prescribed flux} \label{sec:err_est_assess_2}

\begin{figure}[h!]
    \centering
    \begin{subfigure}{0.49\textwidth}
        \includegraphics[trim=3cm 10cm 4cm 9cm, clip=true, width=\textwidth]{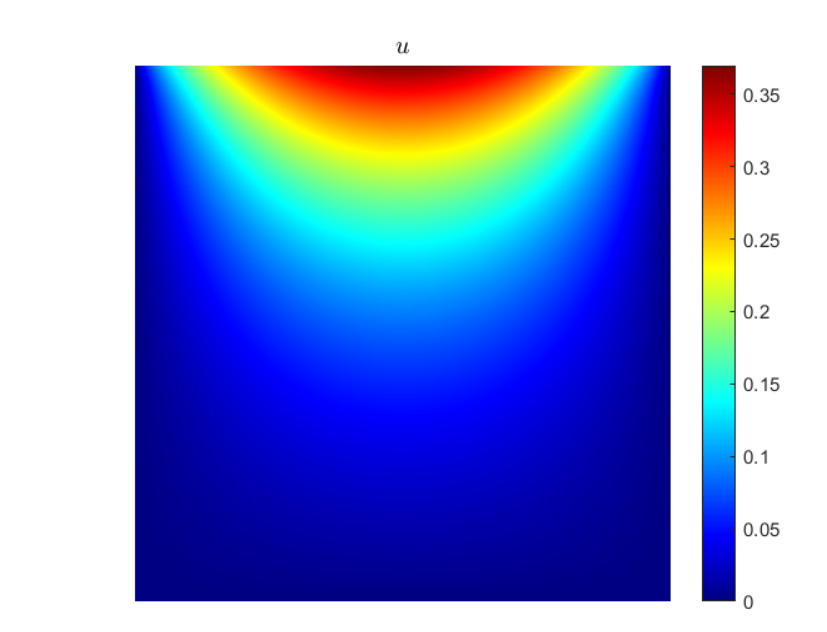}
    \end{subfigure}
    \begin{subfigure}{0.49\textwidth}
        \includegraphics[trim=3cm 10cm 4cm 9cm, clip=true, width=\textwidth]{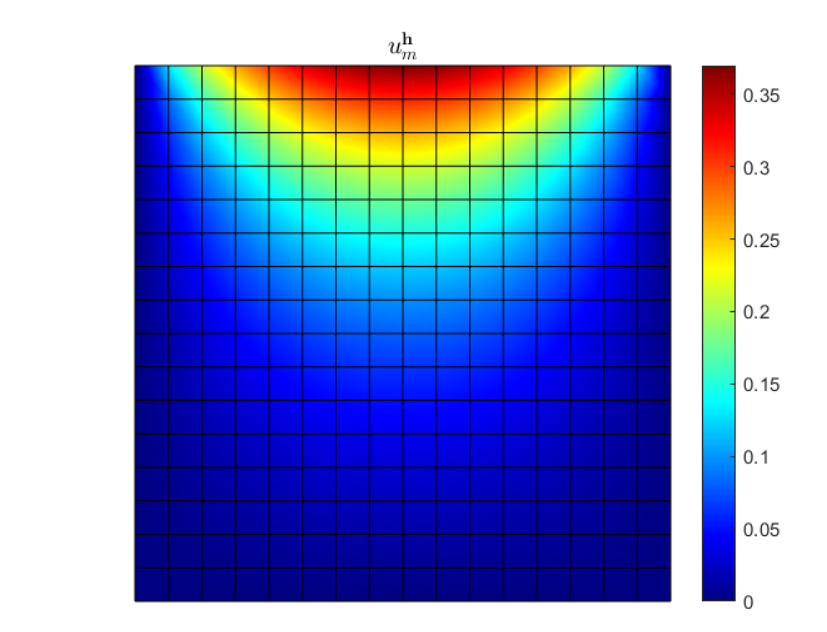}
    \end{subfigure}
    \caption{2D Poisson problem with a prescribed flux: exact solution $u$ (left) and PGD solution with 3 modes (right).}
    \label{fig:sol_2}
\end{figure}

\begin{figure}[h!]
    \centering
    \begin{subfigure}{0.49\textwidth}
        \includegraphics[trim=3cm 10cm 4cm 9cm, clip=true, width=\textwidth]{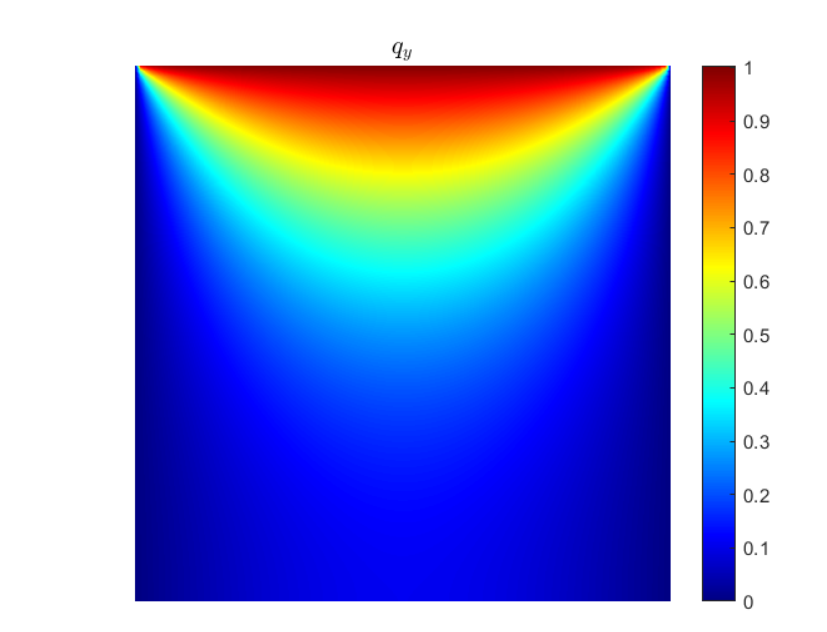}
    \end{subfigure}
    \begin{subfigure}{0.49\textwidth}
        \includegraphics[trim=3cm 10cm 4cm 9cm, clip=true, width=\textwidth]{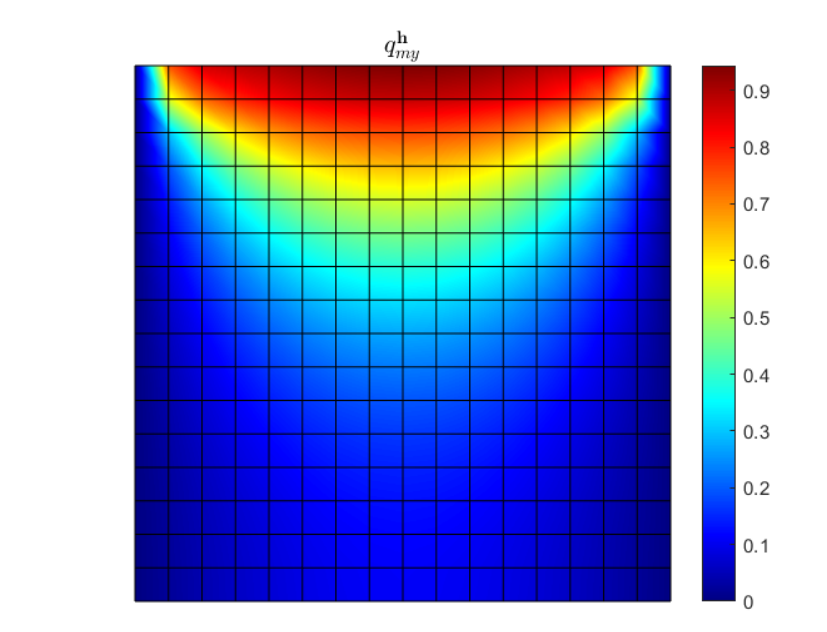}
    \end{subfigure}
    \vfill
    \begin{subfigure}{0.49\textwidth}
        \includegraphics[trim=3cm 10cm 4cm 9cm, clip=true, width=\textwidth]{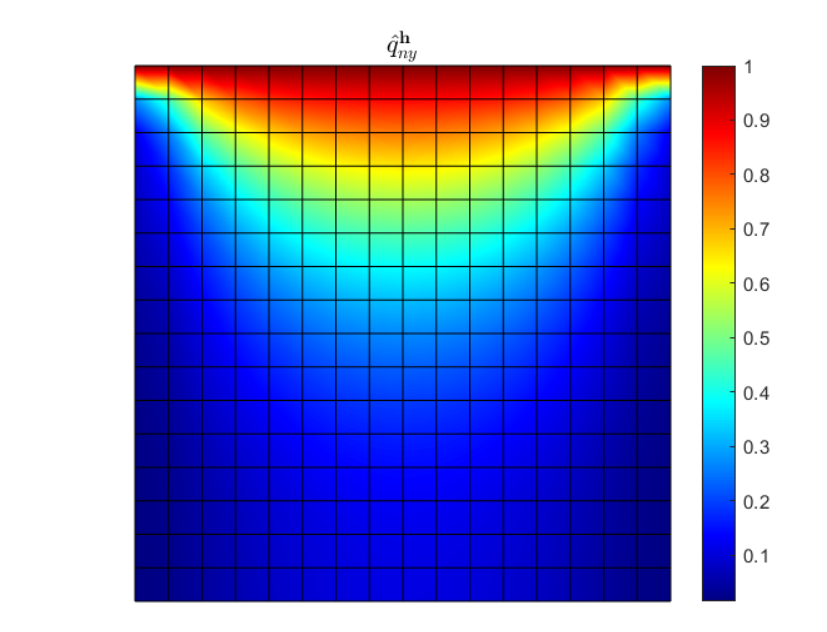}
    \end{subfigure}
    \caption{2D Poisson problem with a prescribed flux: second component of the exact flux $\bs{q}$ (top left), PGD flux $\bs{q}^{\bs{h}}_m$ (top right) and SA PGD flux $\hat{\bs{q}}^{\bs{h}}_n$ with $n=m=3$ modes (bottom).}
    \label{fig:fluxes_y_2}
\end{figure}

We now consider a Poisson problem with a prescribed flux $g=1$ applied on $\Gamma^N=(0,1)\times\{1\}$. Homogeneous Dirichlet boundary conditions are applied on the rest of the boundary and no source term is considered ($f=0$). According to \eqref{eq:q0g}, a flux that equilibrates the external loading is
\begin{equation*}
    \bs{q}_0(x,z) = 
    \begin{pmatrix}
        0 \\
        1
    \end{pmatrix}.
\end{equation*}
The exact solution is also known and reads as
\begin{equation*}
    u(x,z) = \frac{4}{\pi^2} \sum_{\mathrm{odd}\ n=1}^\infty \frac{\sin{(n\pi x)}\sinh{(n\pi z)}}{n^2\cosh{(n\pi)}}.
\end{equation*}
We compute an approximate PGD solution with $N_\omega=N_I=16$ linear elements in each direction,  which is represented in Figure~\ref{fig:sol_2} with 3 modes and compared with the exact solution. The second components of the fluxes $\bs{q}$, $\bs{q}^{\bs{h}}_m$ and $\hat{\bs{q}}^{\bs{h}}_n$ are shown in Figure~\ref{fig:fluxes_y_2}. It should be noted that $\hat{\bs{q}}^{\bs{h}}_n$ is statically admissible. In particular, it satisfies exactly $\hat{\bs{q}}^{\bs{h}}_n\cdot\bs{n}=1$ on $\Gamma^N$, which is not the case for $\bs{q}^{\bs{h}}_m$.

\begin{figure}[h!]
    \centering
    \begin{subfigure}[c]{0.5\textwidth}
        \includegraphics[trim=5cm 9cm 5cm 10cm, clip=true, width=\textwidth]{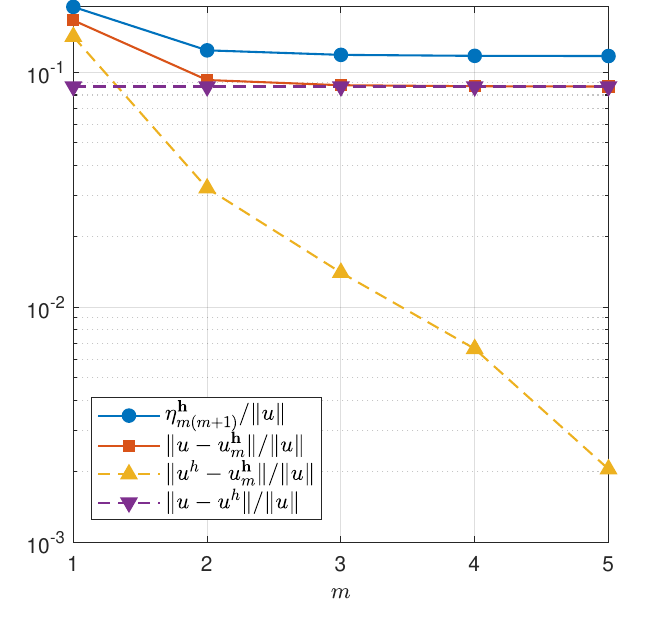}
    \end{subfigure} 
    \hfill
    \begin{subtable}[c]{0.49\textwidth}
        \centering
        \begin{tabular}[b]{c|ccccc}
            \toprule
            $m$ & 1 & 2 & 3 & 4 & 5 \\
            \midrule
            $I_{\text{eff}}$ ($p=1$) & 1.14 & 1.34 & 1.35 & 1.35 & 1.35 \\
            $I_{\text{eff}}$ ($p=2$) & 1.06 & 1.08 & 1.03 & 1.03 & 1.03 \\
            \bottomrule
        \end{tabular}
    \end{subtable}
    \caption{2D Poisson problem with a prescribed flux: exact relative errors and error estimate $\eta_{m\,(m+1)}^{\bs{h}}$ (left),  and associated effectivity index ($I_{\text{eff}}$) with respect to the number of PGD modes and the degree of the elements used to compute $\hat{\bs{q}}_n^{\bs{h}}$ (right).}
    \label{fig:eff_index_2}
\end{figure}

Eventually, the error estimator $\eta_{mn}^{\bs{h}}$ (for $n=m+1$) is shown in Figure~\ref{fig:eff_index_2} with respect to the number $m$ of PGD modes, as well as the energy norm of the exact PGD error and parts of this error due to discretization and reduction, i.e. $\norm{u-u^h}$ and $\norm{u^h-u_m^{\bs{h}}}$ respectively. Here, $u^h$ is the FE solution computed on the 2D mesh $\calT_\omega\otimes\calT_I$ using $Q_1$ elements. The associated effectivity indices are also given in Figure~\ref{fig:eff_index_2}. The error estimator is particularly accurate if quadratic elements are used to compute the modes of $\hat{\bs{q}}^{\bs{h}}_n$. We also take advantage of this test case to show that, as the number of modes increases, the part of the error due to reduction decreases and the overall PGD error reaches a plateau corresponding to the discretization error.

\subsubsection{2D Poisson problem on thin domains} \label{sec:thin_domains}

It should be noted that the PGD technique with separation of space variables is well suited to plate geometries. We are therefore interested in the behaviour of the error estimate in the case where the problem is posed on thin domains. The same Poisson problem as in Section~\ref{sec:err_est_assess_2} is considered, now on the domain $\Omega=(0,1)\times(0,t)$ where the thickness $t$ of the domain may be small ($t\ll 1$).

The effectivity indices associated with $\eta_{mn}^{\bs{h}}$ (for $n=m$) are given in Table~\ref{tab:eff_index_thin} for three different orders of magnitude of $t$. Both PGD solutions were computed with 124 (resp. 16) linear elements in the axial (resp. thickness) direction. It appears that the efficiency index is robust with respect to the thickness $t$.

\begin{table}[h!]
    \centering
    \begin{tabular}{cc|ccc|ccc|cc}
        \toprule
          & & \multicolumn{2}{l}{$t=1$} & & \multicolumn{2}{l}{$t=0.1$} & & \multicolumn{2}{l}{$t=0.01$}\\
        \midrule
        $m$ & & $I_{\text{eff}}$ & $\norm{u-u_m^{\bs{h}}}/\norm{u}$ & & $I_{\text{eff}}$ & $\norm{u-u_m^{\bs{h}}}/\norm{u}$ & & $I_{\text{eff}}$ & $\norm{u-u_m^{\bs{h}}}/\norm{u}$ \\
        \midrule
        1 & & 1.65 &  16~\% & & 1.20 & 8.3~\% & & 1.33 & 3.8~\% \\
        2 & & 1.56 & 6.9~\% & & 1.30 & 2.8~\% & & 1.41 & 3.2~\% \\
        3 & & 1.43 & 6.0~\% & & 1.36 & 2.1~\% & & 1.41 & 3.2~\% \\
        \bottomrule
    \end{tabular}
    \caption{2D Poisson problem on thin domains: influence of slenderness on the efficiency of the error estimator.}
    \label{tab:eff_index_thin}
\end{table}

\section{An adaptive PGD strategy} \label{sec:adaptive_pgd}

In Section~\ref{sec:err_est}, we have developed an error estimator that leads to an upper bound on the overall PGD error. Our goal now is to leverage this estimator to automatically adapt the PGD parameters (the number $m$ of modes and the mesh sizes $h_\omega$ and $h_I$) in order to achieve a prescribed accuracy. 

\subsection{Discussion on possible adaptive approaches}

As already mentioned, two parameters can be adjusted to improve the accuracy of the PGD model:
\begin{itemize}
    \item the number $m$ of modes. If $m$ is too small, the overall error may essentially be a reduction error.
    \item the mesh sizes $h_\omega$ and $h_I$. If the discretization parameters are not chosen correctly, the discretization error can dominate the overall PGD approximation error. 
\end{itemize}
One possible adaptive strategy is to set these parameters in a greedy manner. This approach requires the definition of an error indicator for each error source. To this end, and introducing the FE solution $u^h$ associated with $\calT_\omega \otimes \calT_I$, we can write
\begin{equation} \label{eq:error_decomposition}
    e_{PGD} = u - u_m^{\bs{h}} = \left(u - u^h\right) + \left(u^h - u^{\bs{h}}_m\right) = e^h + e_m,
\end{equation}
where:
\begin{itemize}
    \item $e^h = u - u^h$ is the part of the error due to the discretization;
    \item $e_m = u^h-u^{\bs{h}}_m$ is the part of the error due to reduction, i.e. the number $m$ of modes that is used. 
\end{itemize}
Using the Galerkin orthogonality satisfied by $e^h$ and recalling that $u^{\bs{h}}_m\in V^h$, we have
\begin{equation} \label{eq:error_orthogonality}
    \norm{e_{PGD}}^2 = \norm{e^h}^2 + \norm{e_m}^2.
\end{equation}
Now, assuming that we are able to define an error indicator $\eta_m$ for $e_m$, the discretization error can be estimated by $\eta^{\bs{h}} = \sqrt{\left(\eta_{mn}^{\bs{h}}\right)^2 - \eta_m^2}$. These two error indicators $\eta_m$ and $\eta^{\bs{h}}$ can then be used in a greedy adaptive strategy. In practice, at step $m$, if $\eta_{\bs{h}}<\eta_m$ the mode $m+1$ is computed with the same discretization as for the mode $m$. Otherwise, if $\eta^{\bs{h}}>\eta_m$, the discretization is modified to recompute the mode $m$ with a better accuracy (and the next modes will be computed with this new discretization). 

\begin{remark}
    In order to separate the two error sources, another intermediate solution can be introduced in \eqref{eq:error_decomposition}. We could write 
    \begin{equation*}
        u-u_m^{\bs{h}} = \left(u - u_m\right) + \left(u_m - u_m^{\bs{h}}\right)
    \end{equation*}
    where $u_m$ is the PGD solution without any discretization. Yet another possibility is to write
    \begin{equation*}
        u-u_m^{\bs{h}} = \left(u - u_{m-1}^{\bs{h}} - r_m \otimes s_m\right) + \left(r_m \otimes s_m - r_m^h \otimes s_m^h\right),
    \end{equation*}
    where $u_{m-1}^{\bs{h}}$ is the discrete PGD solution with $m-1$ modes and $(r_m,s_m)$ is solution to \eqref{eq:euler_lagrange}.
    However, in the latter two cases, we do not benefit from the orthogonality property \eqref{eq:error_orthogonality} between the two error sources. This is why we prefer the decomposition \eqref{eq:error_decomposition}. \qed 
\end{remark}

\begin{remark}
    An indicator of the reduction error can be obtained as follows. Noting that $e_m$ satisfies, for any $v\in V^h$,
    \begin{equation*}
        B_1(e_m,v) = R_m^h(v),
    \end{equation*}
    the sum of the modes following the mode $m$ can be seen as an approximation of $e_m$. Considering $M>m$, bounds on the reduction error can be obtained by assuming that $u_{M}^{\bs{h}}$ is a better approximation of $u^h$ than $u_m^{\bs{h}}$. Assuming that there exists a positive constant $C_m<1$ such that,
    \begin{equation} \label{eq:saturation_hyp}
        \norm{e_M} \leq C_m\norm{e_m},
    \end{equation}
    we then have,
    \begin{equation} \label{eq:bounds_red_error}
        \frac{\norm{u^{\bs{h}}_{M}-u^{\bs{h}}_m}}{1+C_m}\leq \norm{e_m} \leq \frac{\norm{u^{\bs{h}}_{M}-u^{\bs{h}}_m}}{1-C_m}.
    \end{equation}
    Inequalities \eqref{eq:bounds_red_error} simply result from \eqref{eq:saturation_hyp} and the triangular inequality. A reduction error indicator $\eta_m$ can then be defined by
    \begin{equation*}
        \eta_m = \norm{u^{\bs{h}}_{M}-u^{\bs{h}}_m}.
    \end{equation*} \qed
\end{remark}

Although greedy, the adaptive strategy described above is not necessarily optimal in terms of computational cost. In fact, in the worst case scenario, the meshes must be modified for each mode computation, which requires recomputing the finite element matrices and projecting the modes already computed onto the new finite element spaces associated to the new meshes. For this reason, in this work, we choose to decouple the two contributions to the overall PGD error. The adaptive strategy adopted here is therefore as follows:
\begin{itemize}
    \item Given a fixed discretization, the approximate PGD solution is enriched with new modes until the reduction error is negligible compared to the discretization error.
    \item If necessary, the discretization is modified and we recompute a PGD solution on the new meshes until the specified tolerance is reached.  
\end{itemize}
This strategy, detailed in the following sections, does not require indicators for each error source. Furthermore, the desired accuracy is expected to be achieved in just a few mesh adaptation steps relative to the number of modes, which reduces the computational costs associated with redefining the finite element matrices. 

\subsection{Adaptivity of PGD modes}

In this section, let us assume that the discretization is fixed. The objective is to define a stopping criterion for the mode enrichment procedure. In view of \eqref{eq:error_decomposition}, for a sufficiently large mode number, the overall PGD error is essentially a discretization error. In other words, as modes are added, the overall PGD error measurement is expected to reach a plateau corresponding to the part of the error due to discretization. This leads us to use our error estimator to define a stopping criterion measuring the stagnation of the PGD error. We enrich the solution by adding new modes as long as
\begin{equation} \label{eq:delta}
    \delta = \left| \frac{\eta_{mn}^{\bs{h}} - \eta_{(m-1)\,n}^{\bs{h}}}{\eta_{(m-1)\,n}^{\bs{h}}} \right| > \varepsilon_{M},
\end{equation}
where $\varepsilon_{M}$ is a user-defined parameter. This criterion measures the relevance of adding a mode with fixed discretization. 

The main stopping criterion of the adaptive process concerns the relative error in energy norm. The objective of the adaptive stategy is to reach the desired accuracy, that is, to guarantee that
\begin{equation*}
    \frac{\norm{u-u_m^{\bs{h}}}}{\norm{u}} < \varepsilon,
\end{equation*}
where $\varepsilon$ is a prescribed tolerance. We indeed prefer a criterion based on relative error rather than a condition based on absolute error, which has less physical meaning. Moreover, the bound $\norm{u-u_m^{\bs{h}}} < \eta_{mn}^{\bs{h}}$ and the triangular inequality imply that
\begin{equation} \label{eq:eta_pgd}
    \frac{\norm{u-u_m^{\bs{h}}}}{\norm{u}} < \eta_{PGD}=\frac{\,\eta_{mn}^{\bs{h}}}{\norm{u_m^{\bs{h}}}-\eta_{mn}^{\bs{h}}},
\end{equation}
we define the main stopping criterion as
\begin{equation} \label{eq:adaptivity_criterion}
    \eta_{PGD} < \varepsilon.
\end{equation}

In the numerical examples below, we set $n=m+1$ and apply the mode adaptation procedure summarized in Algorithm~\ref{alg:mode_adapt}. Note that a maximum number $M$ of PGD modes is introduced to ensure that the procedure has an end. The role of variable \textsf{refineMesh} becomes clear from Algorithm~\ref{alg:adaptive_pgd}.

\begin{algorithm} 
    \SetKwData{RefineMesh}{refineMesh} \SetKwData{Break}{break} \SetKwData{False}{false}
    \KwIn{Two meshes $\calT_\omega$ and $\calT_I$, $\varepsilon$ the tolerance on the error, $\varepsilon_{M}$ the error stagnation rate, $M$ the maximum number of PGD modes.}
    Compute $\hat{\bs{q}}_1^{\bs{h}}$ solving \eqref{eq:dual_pgd_problem} and set $\eta_{01}^{\bs{h}} =\norm{\hat{\bs{q}}_1^{\bs{h}}}_q$\;
    \For{$m=1$ \KwTo $M$}{
    Compute the mode $m$ of $u_m^{\bs{h}}$ solving \eqref{eq:pgd_problem_dis} and get $\norm{u_m^{\bs{h}}}$ \;
    Compute the mode $m+1$ of $\hat{\bs{q}}_{(m+1)}^{\bs{h}}$ solving \eqref{eq:dual_pgd_problem}\;
    Compute the global error estimator $\eta_{m\,(m+1)}^{\bs{h}}=E_{\text{CRE}}\left(u_m^{\bs{h}},\hat{\bs{q}}_{m+1}^{\bs{h}}\right)$ \;
    Deduce $\eta_{PGD}$ from $\eta_{m\,(m+1)}^{\bs{h}}$ and $\norm{u_m^{\bs{h}}}$ as well as $\delta$ from $\eta_{m\,(m+1)}^{\bs{h}}$ and $\eta_{(m-1)\,m}^{\bs{h}}$ (see \eqref{eq:eta_pgd} and \eqref{eq:delta})\;
    \uIf{$\eta_{PGD} < \varepsilon$}{
        Set \RefineMesh = \False \; 
        \Break \;
    }
    \lElseIf{$\displaystyle \delta < \varepsilon_{M}$}{\Break} 
    }
    \caption{Mode adaptation procedure} \label{alg:mode_adapt}
\end{algorithm}

\subsection{Mesh adaptivity and automatic PGD model reduction}

The ultimate goal is not only to determine the number of modes required, but also the finite element partitions $\calT_\omega$ and $\calT_I$, in order to ensure that the relative error in energy norm is below a prescribed tolerance. If the criterion \eqref{eq:adaptivity_criterion} has not been satisfied during the mode adaptation procedure, the discretization needs to be adapted in order to reach the prescribed accuracy. Regarding this task, we follow the standard cycle ESTIMATE, MARK and REFINE. 

As illustrated in Section~\ref{sec:err_est_assess_2}, at the end of the mode adaptation procedure, the overall PGD error is essentially due to discretization. We therefore expect $\eta_{mn}^{\bs{h}}$ to be a good estimate of this discretization error when the PGD solution includes enough modes. We then use local contributions from the error estimate $\eta^{\bs{h}}_{mn}$ to adapt the discretization after the mode adaptation procedure, if necessary.  Since two spatial meshes come into play, we define the following local error estimates, for any element $K$ (resp. $L$) of $\calT_\omega$ (resp. $\calT_I$),
\begin{equation} \label{eq:local_err_ind}
    \eta_K^\omega = \norm{\hat{\bs{q}}_n^{\bs{h}}-\bbA\bs{\nabla}u_m^{\bs{h}}}_{q,\,K\,\times\, I} \quad \eta_L^I = \norm{\hat{\bs{q}}_n^{\bs{h}}-\bbA\bs{\nabla}u_m^{\bs{h}}}_{q,\,\omega\,\times\, L},
\end{equation}
so that we have 
\begin{equation*}
    \left(\eta_{mn}^{\bs{h}}\right)^2 = \sum_{K\,\in\,\calT_\omega} \left(\eta_K^\omega\right)^2 = \sum_{L\,\in\,\calT_I} \left(\eta_L^I\right)^2.
\end{equation*}
Thus, each element $K$ of the mesh $\calT_\omega$ is associated with a local error corresponding to the contribution of the domain $K\times I$ to the total error (and likewise for any element $L$ of $\calT_I$).

In this work, we only consider the $h$-version of mesh adaptation and focus on the local $h$-refinement method. It could also be possible to consider a global $h$-remeshing strategy in the PGD context as shown in \cite{nadal_separated_2015,reis_error_2020}. The elements associated with the largest local error indicators are marked following the maximum strategy (or a variant). We distinguish two cases depending on the dimension of the problem:
\begin{itemize}
    \item For $d=2$, since the quantities $\left(\eta_K^{\omega}\right)_{K\in\calT_\omega}$ and $\left(\eta_L^I\right)_{L\in\calT_I}$ are comparable, we leave ourselves the option of refining in only one direction. Thus, elements $K$ in $\calT_\omega$ such that $\displaystyle \eta_K^{\omega} > \lambda\ \underset{K\in\calT_\omega,\, L\in\calT_I}{\max}\left(\eta_K^{\omega},  \eta_L^I\right)$ for some $\lambda\in (0,1)$ are marked, and likewise for elements $L$ in $\calT_I$;
    \item For $d=3$, the local error indicators  no longer play a similar role. The elements of $\calT_\omega$ and $\calT_I$ are still marked simultaneously but independently. More precisely, elements $K$ in $\calT_\omega$ such that $\displaystyle \eta_K^{\omega} > \lambda_\omega \ \left(\sum_{K\in\calT_\omega}\left(\eta_K^{\omega}\right)^2 /  N_\omega\right)^{1/2}$ are marked, as well as elements $L$ in $\calT_I$ such that $\displaystyle \eta_L^{I} > \lambda_I \ \underset{L\in\calT_I}{\max}\left(\eta_L^I\right)$ where $(\lambda_\omega,\lambda_I)\in (0,1)^2$ are user-chosen parameters.
\end{itemize}

The whole adaptive PGD strategy is presented in Algorithm~\ref{alg:adaptive_pgd}.

\begin{algorithm}
    \SetKwData{RefineMesh}{refineMesh}\SetKwData{True}{true}\SetKwData{False}{false}
    
    Define two initial meshes $\calT_\omega^{0}$ and $\calT_I^0$\;
    Set the refinement parameter(s) $\lambda\in(0,1)$ $\left((\lambda_\omega,\lambda_I)\in(0,1)^2\right)$ and $N$ the maximum number of mesh refinement steps\;
    Set $i = 0$ and \RefineMesh = \True\;
    \While{\RefineMesh}{
        \textbf{SOLVE.} Mode adaptation procedure following Algorithm~\ref{alg:mode_adapt}\;
         Set $i = i+1$ \;
        \If{\RefineMesh}{
            \eIf{$i \leq N$}{
                \textbf{ESTIMATE.} Compute the local PGD error indicators $\left(\eta_K^{\omega}\right)_{K\in\calT_\omega^{i-1}}$ and $\left(\eta_L^I\right)_{L\in\calT_I^{i-1}}$ from \eqref{eq:local_err_ind}\;
                \textbf{MARK.} Mark elements $K$ in $\calT_\omega^{i-1}$ and $L$ in $\calT_I^{i-1}$ associated with largest error indicators\;
                \textbf{REFINE.} Refine marked elements and define the new meshes $\calT_\omega^i$ and $\calT_I^i$\;
            }{
                Set \RefineMesh =\False\;
            }
        }
    }
    \caption{Adaptive PGD procedure} \label{alg:adaptive_pgd}
\end{algorithm}

\subsection{Numerical assessment of the adaptive PGD strategy}

This section is devoted to the evaluation of the adaptive PGD strategy presented above through two test cases. For all the numerical experiments reported here, the following parameters are used: $\epsilon_{FP}=10^{-4}$, $\varepsilon_{M} = 10^{-2}$ and $\lambda=\lambda_\omega=\lambda_I=0.5$. Furthermore, only piecewise linear elements are used for both PGD solutions and $n$ is set to $m+1$. 

\subsubsection{2D Poisson problem with a prescribed flux}

\begin{figure}[h!]
    \centering
    \begin{subfigure}{0.49\textwidth}
        \includegraphics[trim=3cm 9cm 4cm 9cm, clip=true, width=\textwidth]{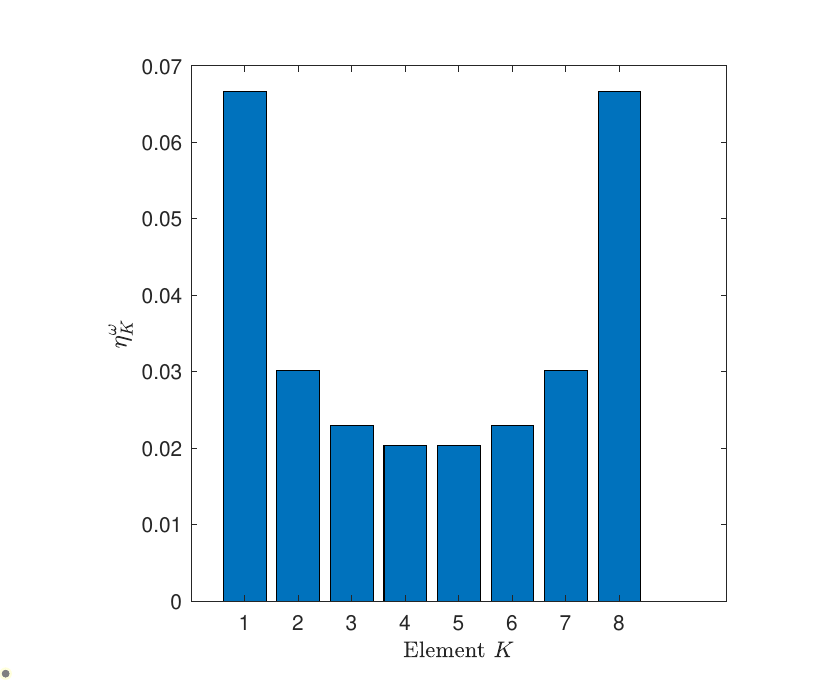}
    \end{subfigure}
    \begin{subfigure}{0.49\textwidth}
        \includegraphics[trim=3cm 9cm 4cm 9cm, clip=true, width=\textwidth]{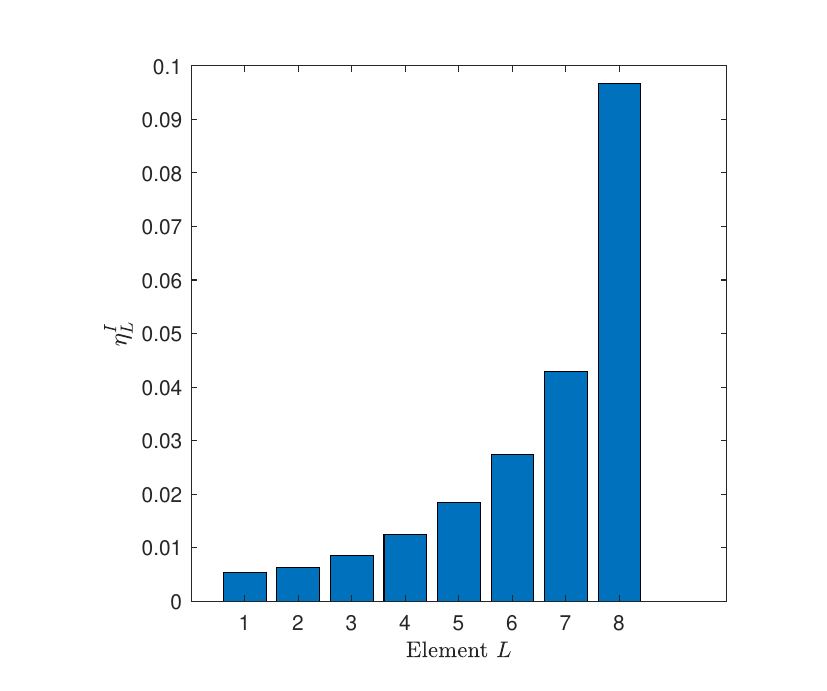}
    \end{subfigure}
    \caption{2D Poisson problem with a prescribed flux: initial values of $\eta_K^\omega$ (left) and $\eta_L^I$ (right).}
    \label{fig:local_errors}
\end{figure}

First, we return to the problem discussed in Section~\ref{sec:err_est_assess_2}. We compute an approximate PGD solution using initial coarse meshes composed of 8 uniform elements. Using the mode adaptation procedure, we obtain a solution comprising 3 modes associated with a relative error estimate $\eta_{PGD}$ of about 28\%. We show the local error estimates in Figure~\ref{fig:local_errors} for this initial uniform discretization. We can see that the error is larger near the boundaries and in particular near the edge where the flux is prescribed. On the other hand, the spatial distribution of the error associated with $\calT_\omega$ exhibits a symmetry, whis is consistent with the symmetry of the problem. 

\begin{figure}[h!]
    \centering
    \begin{subfigure}{0.49\textwidth}
        \includegraphics[trim=3cm 10cm 4cm 9cm, clip=true, width=\textwidth]{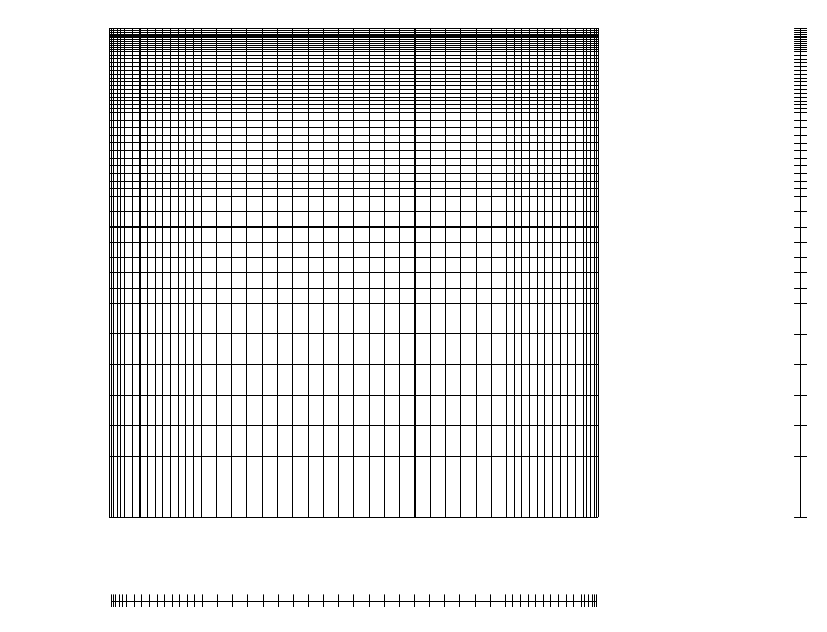}
    \end{subfigure}
    \begin{subfigure}{0.49\textwidth}
        \includegraphics[trim=3cm 9cm 4cm 9cm, clip=true, width=\textwidth]{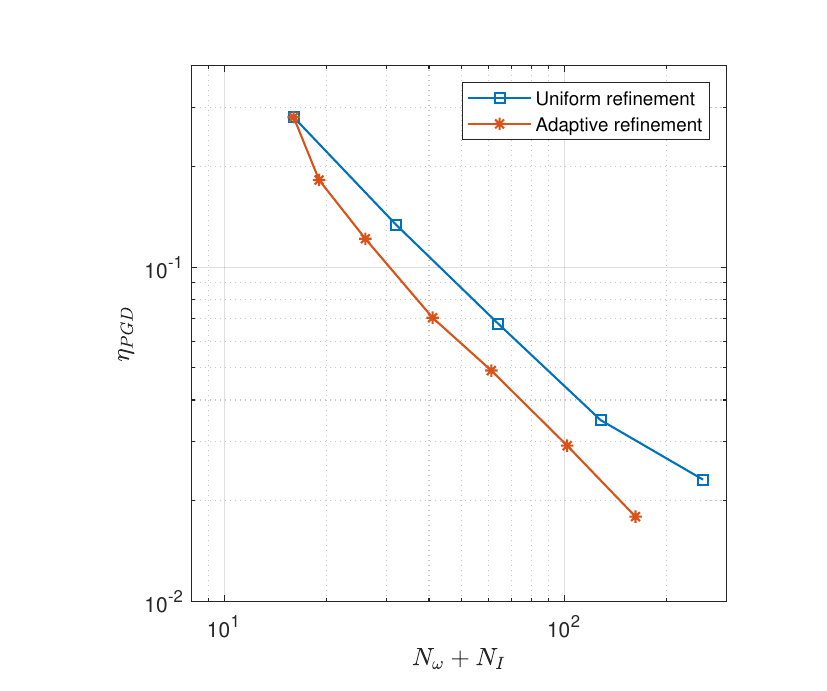}
    \end{subfigure}
    \caption{2D Poisson problem with a prescribed flux: meshes obtained at the end of the adaptive PGD procedure (left) and error convergence for uniform and adaptive mesh refinements (right).}
    \label{fig:refined_meshes}
\end{figure}

We then use the adaptive PGD strategy with a prescribed tolerance of  3\% (i.e. $\varepsilon=0.03$) and from this initial configuration. This results in a PGD solution with 6 modes computed on the meshes shown in Figure~\ref{fig:refined_meshes}. The associated relative error estimate is 2.96\%. As expected, the refinement process leads to finer meshes near the boundaries, particularly near the boundary where the flux is prescribed. Figure~\ref{fig:refined_meshes} also shows the convergence of the relative error estimate $\eta_{PGD}$ with respect to the sum of the degrees of freedom in $\calT_\omega$ and $\calT_I$, for uniform and adaptive refinements. This demonstrates the benefit of using an adaptive refinement algorithm. 

\begin{figure}[h!]
    \centering
        \includegraphics[trim=3cm 9cm 4cm 9cm, clip=true, width=0.5\textwidth]{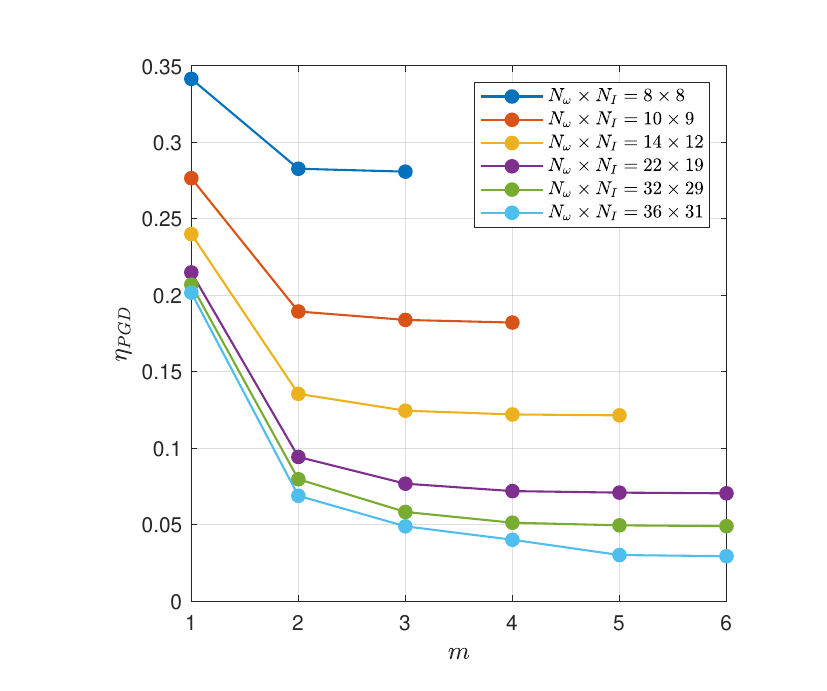}
    \caption{2D Poisson problem with a prescribed flux: history of convergence for the adaptive PGD strategy.}
    \label{fig:conv_hist}
\end{figure}

The history of convergence for the whole adaptive PGD procedure is presented in Figure~\ref{fig:conv_hist}. The evolution of $\eta_{PGD}$ is shown with respect to the number of modes for the different meshes considered during the adaptive procedure. Starting from the initial discretization, the PGD approximation is enriched up to $m=3$ modes. At this stage, the resulting error is mainly due to the spatial discretization. Consequently, the meshes are refined to $10\times 9$ elements, and a new PGD solution is computed with this new discretization. The modal enrichment procedure is restarted from $m=1$ and continued iterativelly to $m=4$. At this stage, the algorithm again triggers a refinement of the meshes. 

As already highlighted, for relatively coarse meshes, the accuracy of the PGD solution is limited by discretization. Thus, the error quickly reaches a plateau when the PGD solution is enriched with new modes. On the other hand, as the meshes become more refined, more modes may be required to reach this plateau.

\subsubsection{2D diffusion problem with discontinuous coefficients}

In this section, we consider \eqref{eq:ref_problem} in the unit square domain $\Omega=(0,1)^2$ with $g=1$ on $(0,1)\times\{1\}$ and $g=0$ on $\{1\}\times (0,1)$. Homogeneous Dirichlet boundary conditions are applied on the rest of the boundary and $f=0$. According to \eqref{eq:q0g},
\begin{equation*}
    \bs{q}_0 = \begin{pmatrix}
        0 \\
        1
    \end{pmatrix}
\end{equation*}
equilibrates the external loading. We consider that $\bbA = \mathbb{I}$ in the first half of the domain ($x<0.5$) and $\bbA = 10\,\mathbb{I}$ in the second half.

\begin{figure}[h!]
    \centering
    \begin{subfigure}{0.49\textwidth}
        \includegraphics[trim=3cm 9cm 4cm 9cm, clip=true, width=\textwidth]{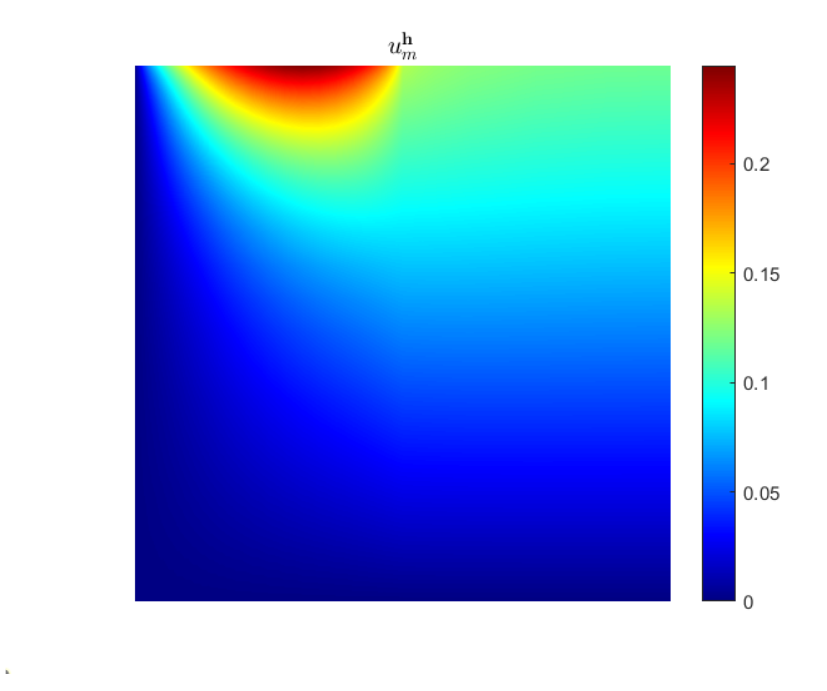}
    \end{subfigure}
    \begin{subfigure}{0.49\textwidth}
        \includegraphics[trim=3cm 9cm 4cm 9cm, clip=true, width=\textwidth]{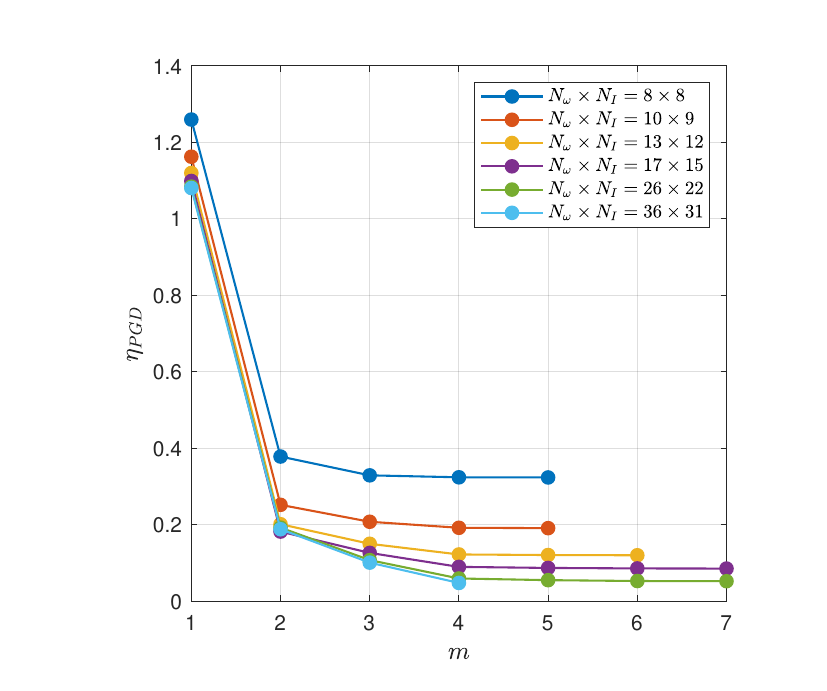}
    \end{subfigure}
    \caption{2D diffusion problem with discontinuous coefficient: PGD solution with 4 modes obtained at the end of the adaptive PGD procedure (left) and history of convergence (right).}
    \label{fig:sol_pgd_3}
\end{figure}

\begin{figure}[h!]
    \centering
    \begin{subfigure}{0.49\textwidth}
        \includegraphics[trim=3cm 10cm 4cm 9cm, clip=true, width=\textwidth]{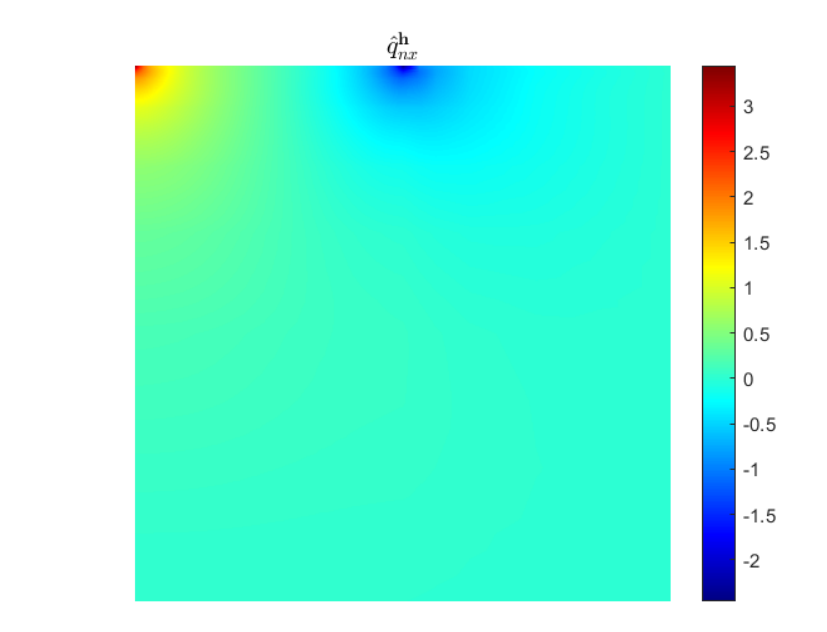}
    \end{subfigure}
    \begin{subfigure}{0.49\textwidth}
        \includegraphics[trim=3cm 10cm 4cm 9cm, clip=true, width=\textwidth]{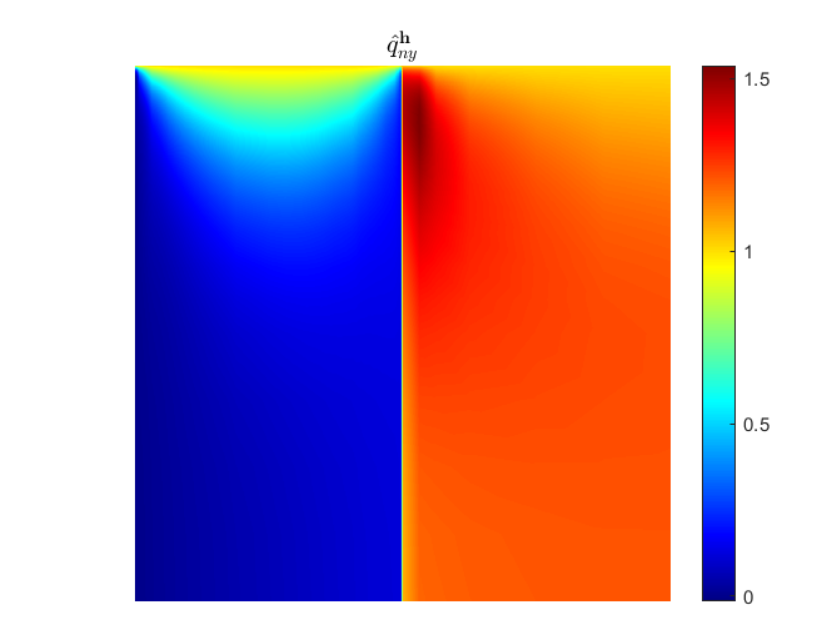}
    \end{subfigure}
    \caption{2D diffusion problem with discontinuous coefficient: components of the SA PGD flux $\hat{\bs{q}}^{\bs{h}}_n$ with 5 modes at the end of the adaptive PGD procedure.}
    \label{fig:flux_pgd_dual_3}
\end{figure}

We use the adaptive PGD algorithm with an initial coarse discretization composed of 8 elements in each direction and with a prescribed error tolerance of 5\%. This results in a PGD solution with 4 modes shown in Figure~\ref{fig:sol_pgd_3}. The value of the relative error estimate $\eta_{PGD}$ decreases from 32\% for the initial configuration to 4.76\%. The history of convergence of the adaptive PGD  procedure is also presented in Figure~\ref{fig:sol_pgd_3}. The final complementary PGD flux $\hat{\bs{q}}^{\bs{h}}_n$ is shown in Figure~\ref{fig:flux_pgd_dual_3}. It should be noted that, in terms of boundary conditions, this flux field is statically admissible. 

\begin{figure}[h!]
    \centering
    \begin{subfigure}{0.49\textwidth}
        \includegraphics[trim=3cm 10cm 4cm 9cm, clip=true, width=\textwidth]{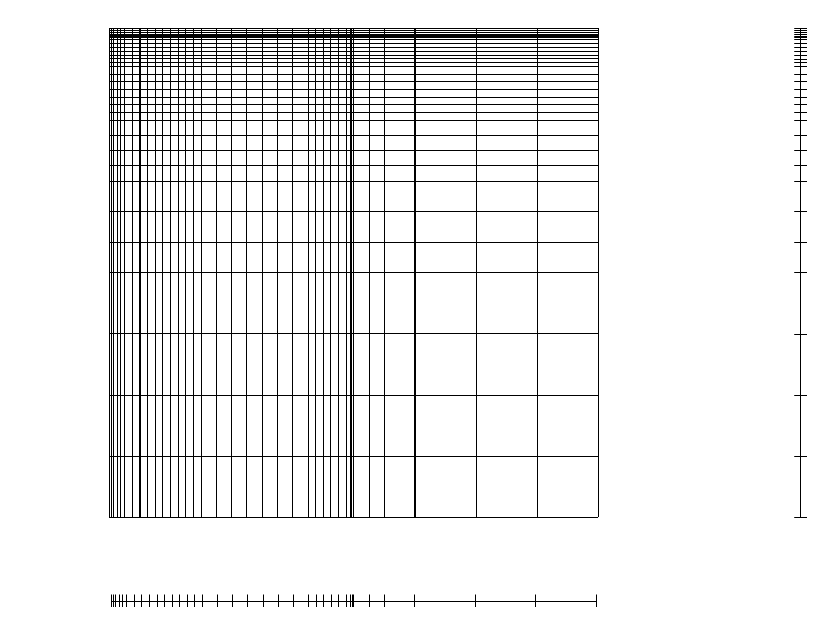}
    \end{subfigure}
    \begin{subfigure}{0.49\textwidth}
        \includegraphics[trim=3cm 9cm 4cm 9cm, clip=true, width=\textwidth]{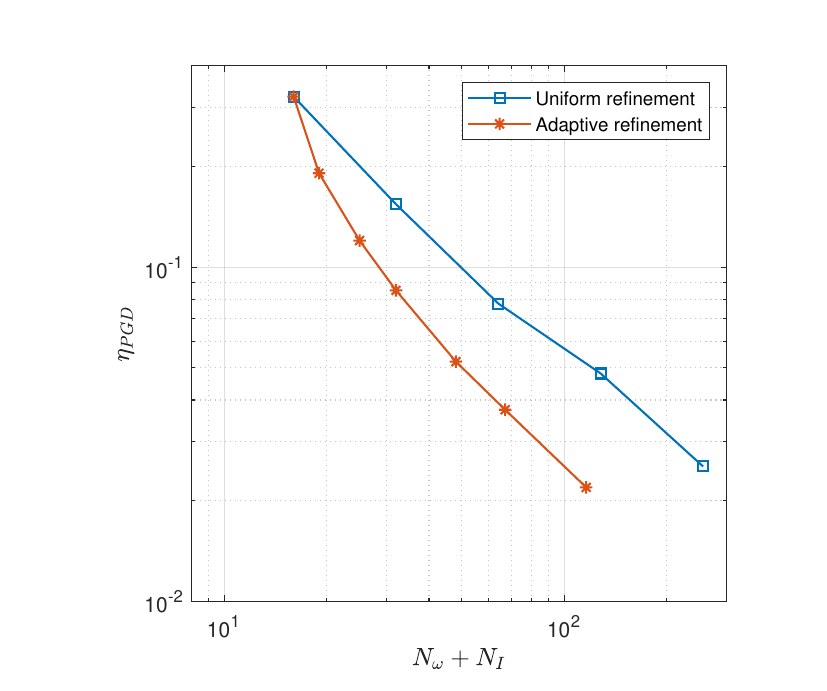}
    \end{subfigure}
    \caption{2D diffusion problem with discontinuous coefficient: meshes obtained at the end of the adaptive PGD procedure (left) and error convergence for uniform and adaptive mesh refinements (right).}
    \label{fig:refined_meshes_2}
\end{figure}

The final meshes (with $N_\omega=36$ and $N_I=31$) are plotted in Figure~\ref{fig:refined_meshes_2} (left). The meshes are particularly refined at the discontinuity and at the left and top boundaries, while being fairly coarse at the bottom right of the domain. Furthermore, compared to uniform refinement, the discretization obtained here is more efficient because it involves fewer degrees of freedom for a given accuracy (Figure~\ref{fig:refined_meshes_2}, right). 

\subsubsection{3D diffusion problem in a laminated plate}

We next consider \eqref{eq:ref_problem} in a plate geometry with $\omega = (-1,1)\times (-1,1)$, $\displaystyle \Omega = \omega\times \left(-\frac{t}{2},\frac{t}{2}\right)$ with $t=0.1$ and $g=1$ on $\displaystyle \omega\times \left\{\frac{t}{2}\right\}$. Homogeneous Dirichlet boundary conditions are applied on the rest of the boundary and $f=0$. By symmetry, the solution is only computed on the top right quarter of the domain. Once again, 
\begin{equation*}
    \bs{q}_0 = \begin{pmatrix}
        \bs{0} \\
        1
    \end{pmatrix}
\end{equation*}
equilibrates the external loading. We also consider that the plate is laminated and composed of three layers, two of them of thickness $t/4$ at the top and bottom of the plate where $\bbA = 0.5\mathbb{I}$ and one in the middle where $\bbA = 5\mathbb{I}$. With regard to the discretization of 2D problems, triangular elements are used. For mesh refinement, each marked triangle is divided into four triangles. Unmarked triangles can also be refined to avoid hanging nodes. 

\begin{figure}[h!]
    \centering
    \begin{subfigure}{0.49\textwidth}
        \includegraphics[trim=3cm 9cm 4cm 9cm, clip=true, width=\textwidth]{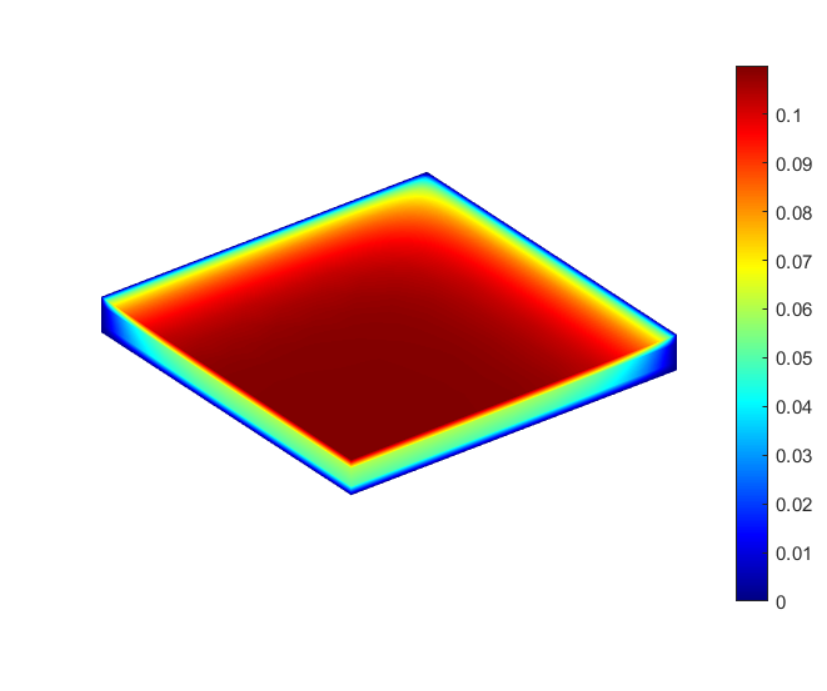}
    \end{subfigure}
    \begin{subfigure}{0.49\textwidth}
        \includegraphics[trim=3cm 9cm 4cm 9cm, clip=true, width=\textwidth]{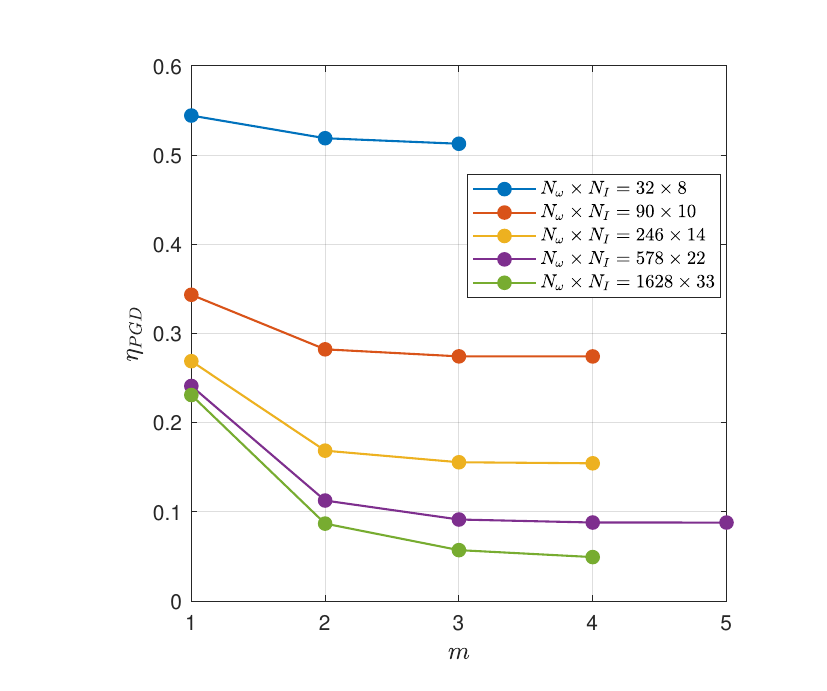}
    \end{subfigure}
    \caption{3D diffusion problem in a laminated plate: PGD solution with 3 modes obtained at the end of the adaptive PGD procedure (left) and history of convergence (right).}
    \label{fig:sol_pgd_4}
\end{figure}

\begin{figure}[h!]
    \centering
    \begin{subfigure}{0.49\textwidth}
        \includegraphics[trim=3cm 10cm 4cm 9cm, clip=true, width=\textwidth]{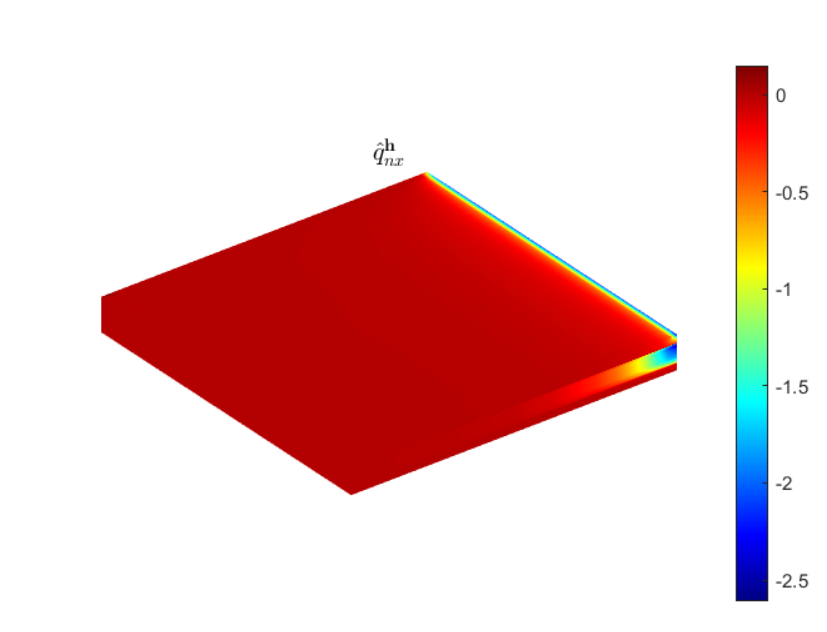}
    \end{subfigure}
    \begin{subfigure}{0.49\textwidth}
        \includegraphics[trim=3cm 10cm 4cm 9cm, clip=true, width=\textwidth]{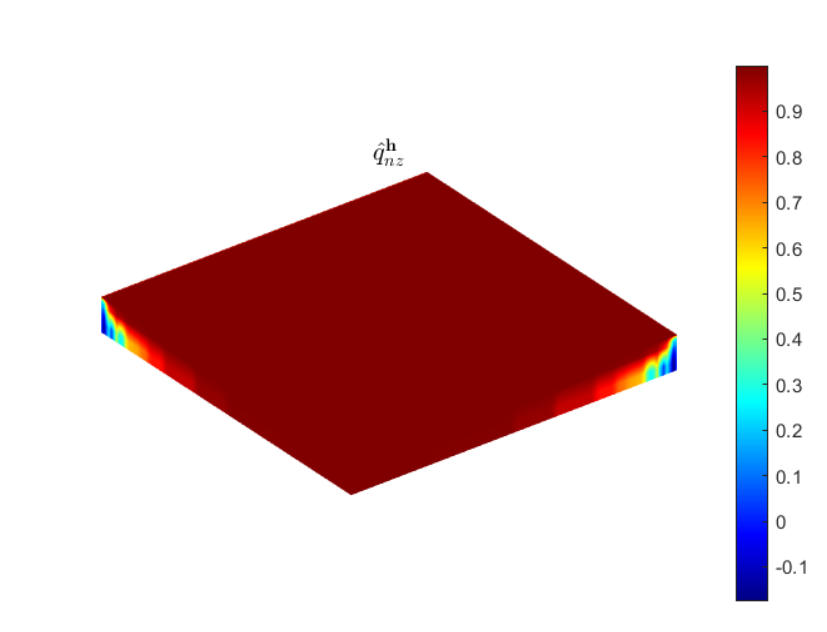}
    \end{subfigure}
    \caption{3D diffusion problem in a laminated plate: components $x$ and $y$ of the SA PGD flux $\hat{\bs{q}}^{\bs{h}}_n$ with 4 modes at the end of the adaptive PGD procedure.}
    \label{fig:flux_pgd_dual_4}
\end{figure}

We use the adaptive PGD algorithm starting with a initial coarse discretization and a prescribed error tolerance of 5\%. This results in a PGD solution with 3 modes shown in Figure~\ref{fig:sol_pgd_4} (left). The value of the relative error estimate $\eta_{PGD}$ decreases from 51\% to 4.94\%. The history of convergence of the adaptive PGD procedure is presented in Figure~\ref{fig:sol_pgd_4} (right). Two components of the final SA PGD flux are shown in Figure~\ref{fig:flux_pgd_dual_4}, where it can be verified that it exactly satisfies the boundary conditions. 

\begin{figure}[h!]
    \centering
    \begin{subfigure}{0.49\textwidth}
        \includegraphics[trim=3cm 9cm 4cm 9cm, clip=true, width=\textwidth]{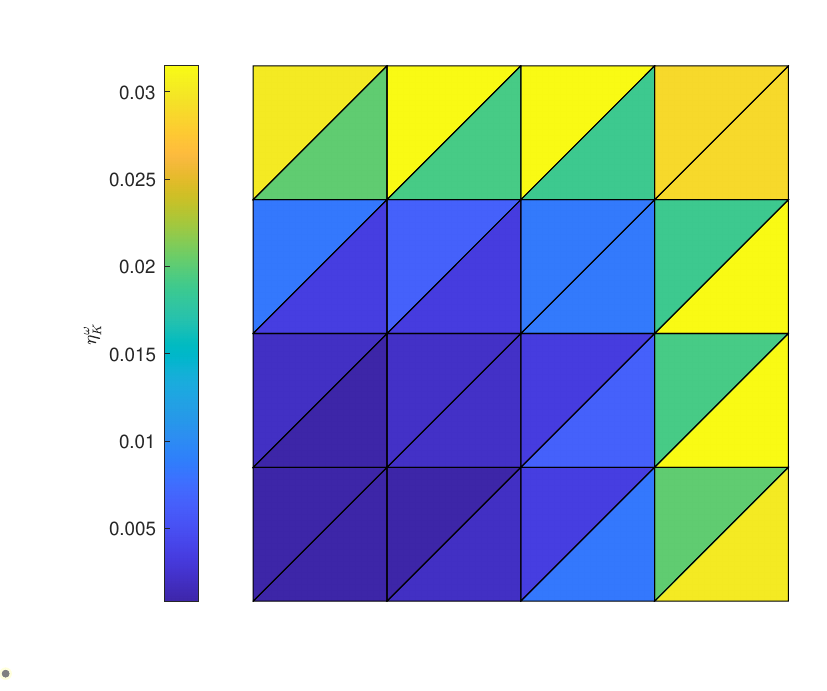}
    \end{subfigure}
    \begin{subfigure}{0.49\textwidth}
        \includegraphics[trim=3cm 9cm 4cm 9cm, clip=true, width=\textwidth]{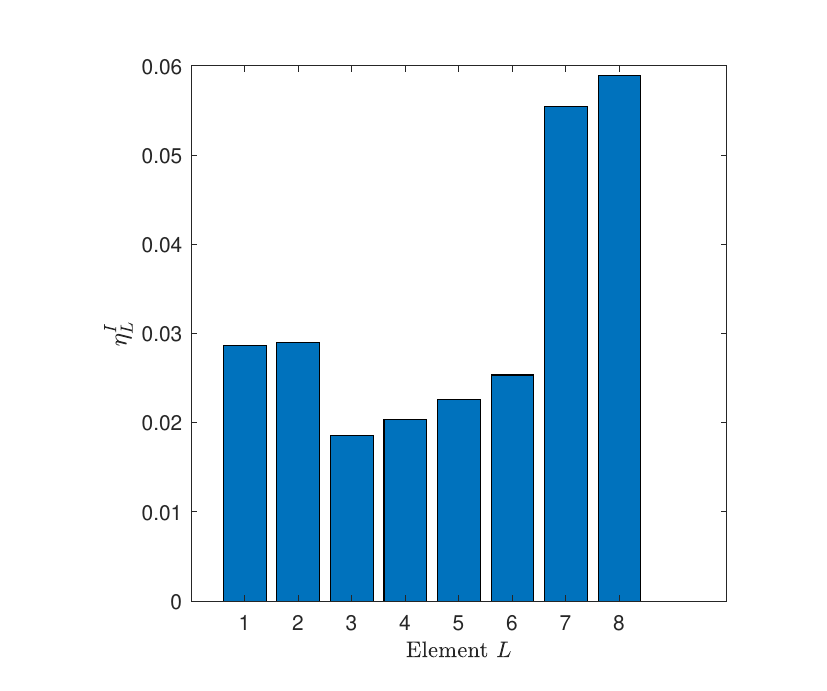}
    \end{subfigure}
    \caption{3D diffusion problem in a laminated plate: initial values of $\eta_K^\omega$ (left) and $\eta_L^I$ (right) when considering $m=3$ modes.}
    \label{fig:local_errors_2}
\end{figure}

\begin{figure}[h!]
    \centering
    \begin{subfigure}{0.49\textwidth}
        \includegraphics[trim=3cm 9cm 4cm 9cm, clip=true, width=\textwidth]{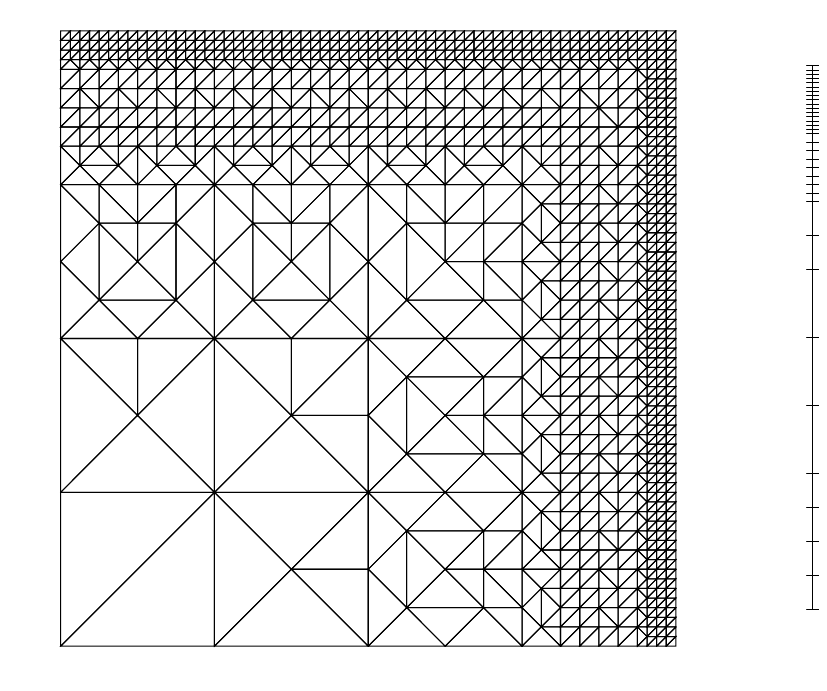}
    \end{subfigure}
    \begin{subfigure}{0.49\textwidth}
        \includegraphics[trim=3cm 9cm 4cm 9cm, clip=true, width=\textwidth]{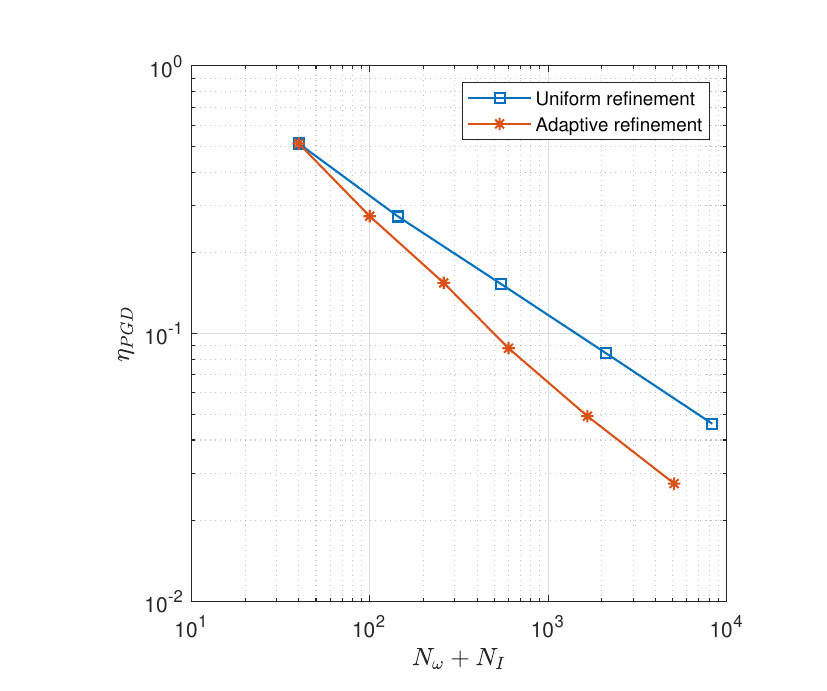}
    \end{subfigure}
    \caption{3D diffusion problem in a laminated plate: meshes obtained at the end of the adaptive PGD procedure (left) and error convergence for uniform and adaptive mesh refinement (right).}
    \label{fig:refined_meshes_3}
\end{figure}

As shown in Figure~\ref{fig:local_errors_2}, when using the initial meshes the error is larger near the Dirichlet boundaries and where the non-zero flux is applied. The adaptive procedure then leads to refining the meshes in these areas. The final meshes (with $N_\omega = 1628$ and $N_I = 33$) are plotted in Figure~\ref{fig:refined_meshes_3} (left). 

\section{Goal-oriented error estimation for the PGD approximation} \label{sec:goal_oriented}

In this section, we show how to extend the previous tools to the framework of goal-oriented error estimation. It is indeed well known that controlling and adapting a numerical solution by measuring the error in energy norm is not optimal when one is interested in a specific quantity of interest.  
 
We consider a quantity of interest $Q(u)$ defined from a linear and continuous functional $Q:V \rightarrow \mathbb{R}$. This functional $Q$ is defined globally using extraction functions $\td{f}_Q$, $\td{\bs{q}}_Q$, $\td{g}_Q$ and $\td{u}_Q$:
\begin{equation}
    Q(v) = \int_\Omega \left(\td{f}_Q \,v + \td{\bs{q}}_Q\cdot\bs{\nabla}v \right) + \int_{\Gamma^N} \td{g}_Q\,v + \int_\Omega \bbA \bs{\nabla}\td{u}_Q\cdot \bs{\nabla} v.
\end{equation}
The objective of the following is therefore to estimate the PGD error in quantity of interest $Q\left(e_{PGD}\right)=Q(u)-Q\left(u_m^{\bs{h}}\right)$. 

\subsection{Adjoint problem and upper bound on the error in quantity of interest}

To this end, the adjoint problem associated with the functional $Q$ is introduced \cite{becker_optimal_2001}. It consists in finding $\tilde{u}\in V$ such that,
\begin{equation*}
    \forall\, v\in V, \quad B_1^*(\td{u},v) = Q(v),
\end{equation*}
where $B_1^*$ is the adjoint operator of $B_1$ defined by $B_1^*(u,v)=B_1(v,u)$. Since $B_1$ is symmetric, we actually have $B_1^*=B_1$.

It is then possible to compute an approximate PGD solution $\td{u}_m^{\bs{h}}$ to the adjoint problem. 

\begin{remark}
    For the sake of clarity, we assume that the solution to the adjoint problem is approximated using the same discretization and the same number of modes as the primal solution. However, it could be possible to relax these constraints. \qed
\end{remark}

We are now able to bound the error on $Q(u)$. Since $e_{PGD}\in V$, we have
\begin{equation} \label{eq:err_bound_calc}
    Q(e_{PGD}) = B_1^*(\td{u},e_{PGD}) = B_1(e_{PGD},\td{u}) = B_1\left(e_{PGD},\td{u}-\td{u}^{\bs{h}}_m\right) + B_1\left(e_{PGD},\td{u}^{\bs{h}}_m\right) =B_1(e_{PGD},\td{e}_{PGD}) + B_1\left(e_{PGD},\td{u}^{\bs{h}}_m\right),
\end{equation}
where $\td{e}_{PGD}=\td{u}-\td{u}^{\bs{h}}_m$ is the PGD error for the adjoint solution. Furthermore, the last term in \eqref{eq:err_bound_calc} is written as
\begin{equation} \label{eq:corr_term}
    B_1\left(e_{PGD},\td{u}^{\bs{h}}_m\right)=L\left(\td{u}^{\bs{h}}_m\right)-B_1\left(u^{\bs{h}}_m,\td{u}^{\bs{h}}_m\right)=R_m\left(\td{u}^{\bs{h}}_m\right)
\end{equation} 
and is a computable term. Thus, by using \eqref{eq:corr_term} and the Cauchy-Schwarz inequality, we deduce from \eqref{eq:err_bound_calc} that 
\begin{equation} \label{eq:err_bound_1}
    \left|Q(e_{PGD})-R_m\left(\td{u}^{\bs{h}}_m\right)\right| \leq \norm{e_{PGD}} \norm{\td{e}_{PGD}},
\end{equation}
where $R_m(\td{u}^{\bs{h}}_m)$ can be seen as a correction term for the quantity of interest. 

\subsection{Goal-oriented PGD error estimator and adaptive strategy}

An estimator of the PGD error in quantity of interest follows directly from \eqref{eq:err_bound_1} using the tools presented in Section~\ref{sec:err_est} \cite{ladeveze_strict_2008,chamoin_non-intrusive_2008}. Introducing the space of equilibrated fluxes for the adjoint problem
\begin{equation*}
    \widetilde{W} = \left\{\td{\bs{p}}\in \left[L^2(\Omega)\right]^d,\ \forall v\in V, \ \int_\Omega \td{\bs{p}}\cdot \bs{\nabla} v = Q(v) \right\},
\end{equation*}
a flux $\hat{\td{\bs{q}}}_n^{\bs{h}} \in\widetilde{W}$ can be computed using the complementary PGD approach detailed in Section~\ref{sec:complementary_PGD}. According to \eqref{eq:err_bound_1} and the property of the CRE functional, we now have
\begin{equation} \label{eq:goal_err_est_1}
    \left|Q(e_{PGD})-R_m\left(\td{u}^{\bs{h}}_m\right)\right| \leq E_{\text{CRE}}\left(u_m^{\bs{h}},\hat{\bs{q}}_n^{\bs{h}}\right) E_{\text{CRE}}\left(\td{u}_m^{\bs{h}},\hat{\td{\bs{q}}}_n^{\bs{h}}\right). 
\end{equation}
Note that a more accurate error bound than \eqref{eq:goal_err_est_1} can be obtained by starting from \eqref{eq:err_bound_calc}. Using the fact that $\hat{\td{\bs{q}}}^{\bs{h}}_n \in \widetilde{W}$ and $u-u^{\bs{h}}_m\in V$, we write:
\begin{equation*}
\begin{aligned}
    Q(e_{PGD}) - R_m\left(\td{u}^{\bs{h}}_m\right) = B_1\left(u-u^{\bs{h}}_m,\td{u}-\td{u}^{\bs{h}}_m\right) = \int_\Omega \nabla\left(u-u^{\bs{h}}_m\right)\cdot \left(\td{\bs{q}}-\bbA\nabla \td{u}^{\bs{h}}_m\right) = \int_\Omega \nabla\left(u-u^{\bs{h}}_m\right)\cdot \left(\hat{\td{\bs{q}}}^{\bs{h}}_n-\bbA\nabla \td{u}^{\bs{h}}_m\right).
\end{aligned}
\end{equation*}
We thus have
\begin{equation} \label{eq:err_bound_calc_bis}
    \begin{aligned}
       Q(e_{PGD}) - R_m\left(\td{u}^{\bs{h}}_m\right) &= \int_\Omega \bbA^{-1}\left(\bs{q}-\bs{q}^{\bs{h}}_m\right)\cdot \left(\hat{\td{\bs{q}}}^{\bs{h}}_n-\bbA\nabla \td{u}^{\bs{h}}_m\right) \\
       & = \int_\Omega \bbA^{-1}(\bs{q}-\bs{q}^*)\cdot \left(\hat{\td{\bs{q}}}^{\bs{h}}_n-\bbA\nabla \td{u}^{\bs{h}}_m\right) + C^{\bs{h}}_m,
    \end{aligned}
\end{equation}
where $\displaystyle \bs{q}^* = \frac{1}{2}\left(\hat{\bs{q}}^{\bs{h}}_n + \bs{q}^{\bs{h}}_m\right)$ and $\displaystyle C^{\bs{h}}_m = \frac{1}{2}\int_\Omega \bbA^{-1}\left(\hat{\bs{q}}^{\bs{h}}_n-\bs{q}^{\bs{h}}_m\right)\cdot \left(\hat{\td{\bs{q}}}^{\bs{h}}_n-\bbA\nabla \td{u}^{\bs{h}}_m\right)$. By applying the Cauchy-Schwarz inequality and using the hypercircle property \eqref{eq:hypercircle}, we obtain from \eqref{eq:err_bound_calc_bis} that
\begin{equation} \label{eq:goal_err_est_2}
    \left|Q(e_{PGD})- \overline{C}^{\bs{h}}_m \right| \leq \frac{1}{2} E_{\text{CRE}}\left(u_m^{\bs{h}},\hat{\bs{q}}_n^{\bs{h}}\right) E_{\text{CRE}}\left(\td{u}_m^{\bs{h}},\hat{\td{\bs{q}}}_n^{\bs{h}}\right),
\end{equation}
where $\displaystyle \overline{C}^{\bs{h}}_m = R_m\left(\td{u}^{\bs{h}}_m\right) + C^{\bs{h}}_m = \frac{1}{2} \int_\Omega \bbA^{-1}\left(\hat{\bs{q}}^{\bs{h}}_n-\bs{q}^{\bs{h}}_m\right)\cdot \left(\hat{\td{\bs{q}}}^{\bs{h}}_n+\bbA\nabla \td{u}^{\bs{h}}_m\right)$. 

In view of \eqref{eq:goal_err_est_2}, we define our goal-oriented PGD error estimator $\eta^Q_{PGD}$ by (see \cite{chamoin_introductory_2022})
\begin{equation}
   \eta^Q_{PGD} = \max_{\theta = \pm 1} \left|\overline{C}^{\bs{h}}_m + \frac{\theta}{2} E_{\text{CRE}}\left(u_m^{\bs{h}},\hat{\bs{q}}_n^{\bs{h}}\right) E_{\text{CRE}}\left(\td{u}_m^{\bs{h}},\hat{\td{\bs{q}}}_n^{\bs{h}}\right) \right|
\end{equation}
and have $ Q(e_{PGD}) \leq \eta^Q_{PGD}$. 
Based on this estimate, an adaptive PGD strategy similar to the one presented in Section~\ref{sec:adaptive_pgd} is used to achieve the desired accuracy on the quantity of interest. For a given discretization, the PGD solutions are enriched with new modes until the error estimator $\eta^Q_{PGD}$ stagnates. If necessary, the meshes are adapted based on the local contributions of $\eta^Q_{PGD}$. Denoting $\theta_{max}$ the maximizer in \eqref{eq:goal_err_est_2}, these local contributions are in the form, for any element $K$ (resp. $L$) of $\calT_\omega$ (resp. $\calT_I$), of
\begin{equation*}
    \overline{C}^{\bs{h}}_{m,\, K\times I} + \frac{\theta_{max}}{2}\, \eta_{K\times I}^Q \quad \text{and} \quad  \overline{C}^{\bs{h}}_{m,\,\omega\times L} + \frac{\theta_{max}}{2}\, \eta_{\omega\times L}^Q,
\end{equation*}
where $\displaystyle \overline{C}^{\bs{h}}_{m,\,\square} = \frac{1}{2} \int_{\square} \bbA^{-1}\left(\hat{\bs{q}}^{\bs{h}}_n-\bs{q}^{\bs{h}}_m\right)\cdot \left(\hat{\td{\bs{q}}}^{\bs{h}}_n+\bbA\nabla \td{u}^{\bs{h}}_m\right)$ and
\begin{equation*}
    \eta_\square^Q = \sqrt{\frac{1}{2}E_{\text{CRE}}^2\left(u_m^{\bs{h}},\hat{\bs{q}}_n^{\bs{h}}\right)\norm{\hat{\td{\bs{q}}}_n^{\bs{h}}-\bbA\bs{\nabla}\td{u}_m^{\bs{h}}}_{q,\,\square}^2 + \frac{1}{2}\norm{\hat{\bs{q}}_n^{\bs{h}}-\bbA\bs{\nabla}u_m^{\bs{h}}}_{q,\,\square}^2 E_{\text{CRE}}^2\left(\td{u}_m^{\bs{h}},\hat{\td{\bs{q}}}_n^{\bs{h}}\right)}.
\end{equation*}
These local contribution are thus defined in such a way that
\begin{equation*}
    \eta^Q_{PGD} = \left| \sum_{K\in \calT_\omega} \overline{C}^{\bs{h}}_{m,\, K\times I} + \frac{\theta_{max}}{2} \sqrt{\sum_{K\in \calT_\omega}\left(\eta_{K\times I}^Q\right)^2} \right| =  \left| \sum_{L\in \calT_I} \overline{C}^{\bs{h}}_{m,\, \omega\times L} + \frac{\theta_{max}}{2} \sqrt{\sum_{L\in \calT_I}\left(\eta_{\omega\times L}^Q\right)^2} \right |.
\end{equation*}

\subsubsection{Numerical illustration}

To illustrate the above goal-oriented PGD strategy, we consider a problem already addressed in \cite{giraldi_tensor_2013} with tensor-based methods. This is a boundary value problem, which corresponds to the diffusion problem \eqref{eq:ref_problem} with $\Gamma^N=\partial\Omega$, $f=0$ and $g=\bs{e}_\alpha\cdot\bs{n}$ where $1\leq\alpha\leq d$. Introducing the space $H^1_*(\Omega)$ of zero-mean $H^1(\Omega)$ functions, the weak form of the problem consists in finding $u\in H^1_*(\Omega)$ such that,
\begin{equation*}
    \forall\, v\in H^1_*(\Omega), \quad \int_\Omega \bbA \bs{\nabla}u\cdot \bs{\nabla} v = \int_{\partial\Omega} \left(\bs{e}_\alpha\cdot \bs{n}\right) v.
\end{equation*}

\begin{remark}
    Such a problem appears in homogenization, and it actually corresponds to the corrector problem (the domain $\Omega$ is then a representative volume element of a heterogeneous material) when considering Neumann boundary conditions (see e.g. \cite[Section 3.1]{lebris_special_2016} and \cite[Section 1.3.3]{bourgeat_approximations_2004}). \qed
\end{remark}

This problem is reformulated as follows: find $u\in H^1(\Omega)$ such that,
\begin{equation} \label{eq:hom_problem}
    \forall\, v\in H^1(\Omega), \quad \int_\Omega \bbA \bs{\nabla}u\cdot \bs{\nabla} v  + \gamma \int_\Omega u \int_\Omega v = \int_{\partial\Omega} \left(\bs{e}_\alpha \cdot \bs{n}\right) v,
\end{equation}
where $\gamma>0$. In this context, the quantities of interest are 
\begin{equation*}
    Q_\alpha(u) = \frac{1}{|\Omega|} \int_\Omega \bs{\nabla} u \cdot \bs{e}_\alpha,
\end{equation*}
from which the effective (homogenized) properties of the material can be extracted (see \cite{giraldi_tensor_2013} for details). Note that by integrating by parts, we obtain that
\begin{equation*}
    Q_\alpha(u) = \frac{1}{|\Omega|} \int_{\partial\Omega} \left(\bs{e}_\alpha \cdot \bs{n}\right) u,
\end{equation*}
from which we deduce that $\displaystyle \td{u} = \frac{1}{|\Omega|}u$. Thus, the adjoint problem does not need to be solved here. 

Inspired by this work, we set $\Omega = (0,1)^2$ and consider the problem \eqref{eq:hom_problem} with $\bbA(x,z) = A(x,z)\mathbb{I}$, where $A(x,z) = 10\, \chi_{[0,\frac{1}{4}]\cup [\frac{3}{4},1]}(x) + \chi_{[\frac{1}{4},\frac{3}{4}]}(x)$. 

\begin{figure}[h!]
    \centering
    \begin{subfigure}{0.49\textwidth}
        \includegraphics[trim=3cm 9cm 4cm 9cm, clip=true, width=\textwidth]{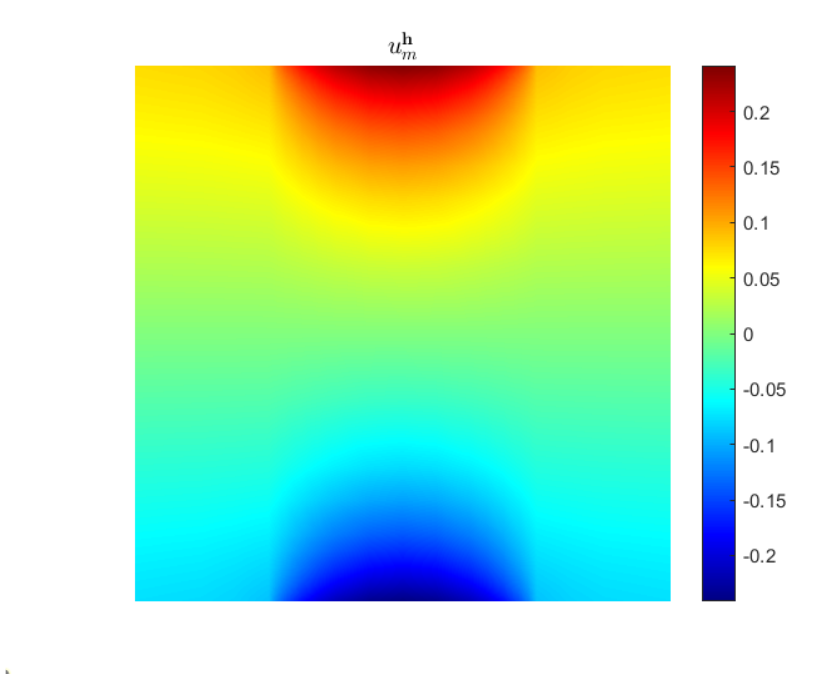}
    \end{subfigure}
    \begin{subfigure}{0.49\textwidth}
        \includegraphics[trim=3cm 10cm 4cm 9cm, clip=true, width=\textwidth]{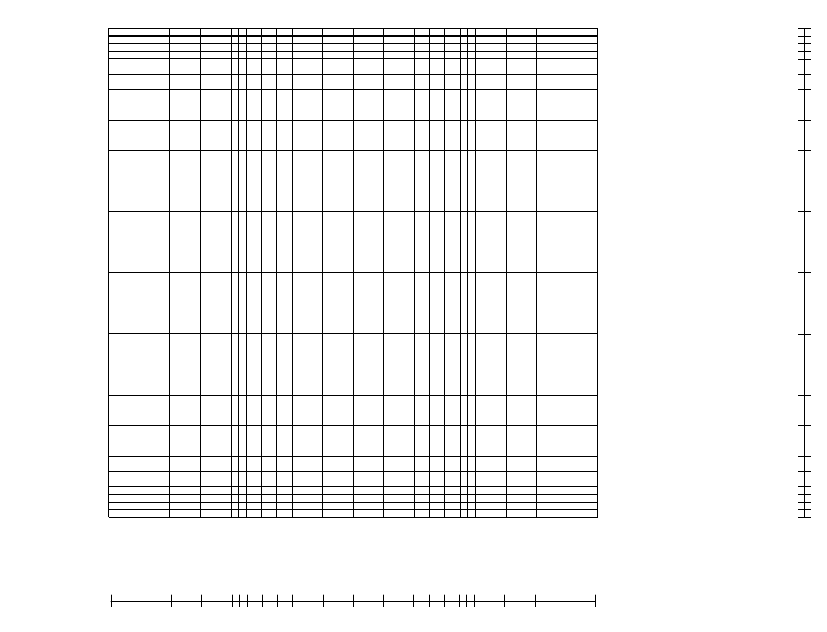}
    \end{subfigure}
    \caption{2D homogenization problem: PGD solution $u^{\bs{h}}_m$ (left) and meshes (right) obtained at the end of the goal-oriented adaptive PGD procedure.}
    \label{fig:sol_pgd_hom}
\end{figure}

\begin{figure}[h!]
    \centering
    \begin{subfigure}{0.49\textwidth}
        \includegraphics[trim=3cm 10cm 4cm 9cm, clip=true, width=\textwidth]{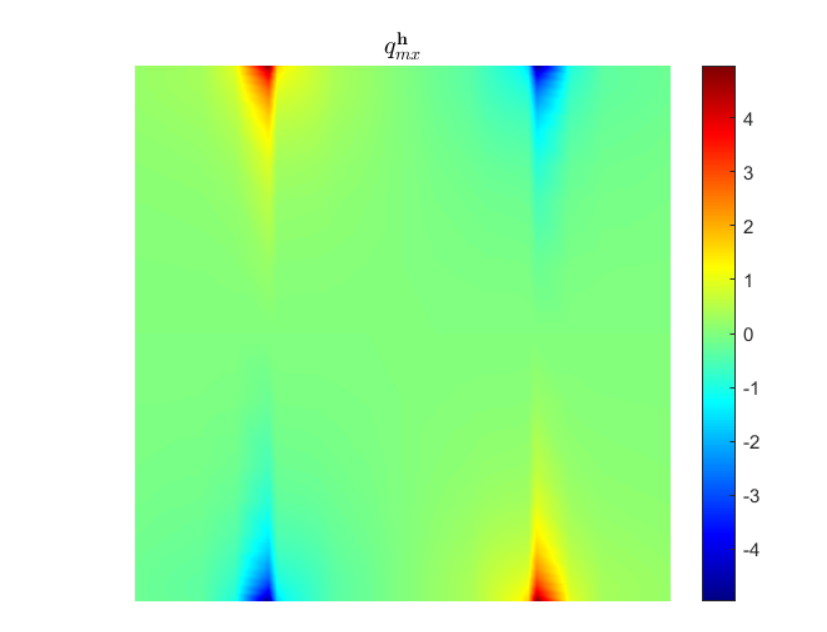}
    \end{subfigure}
    \begin{subfigure}{0.49\textwidth}
        \includegraphics[trim=3cm 10cm 4cm 9cm, clip=true, width=\textwidth]{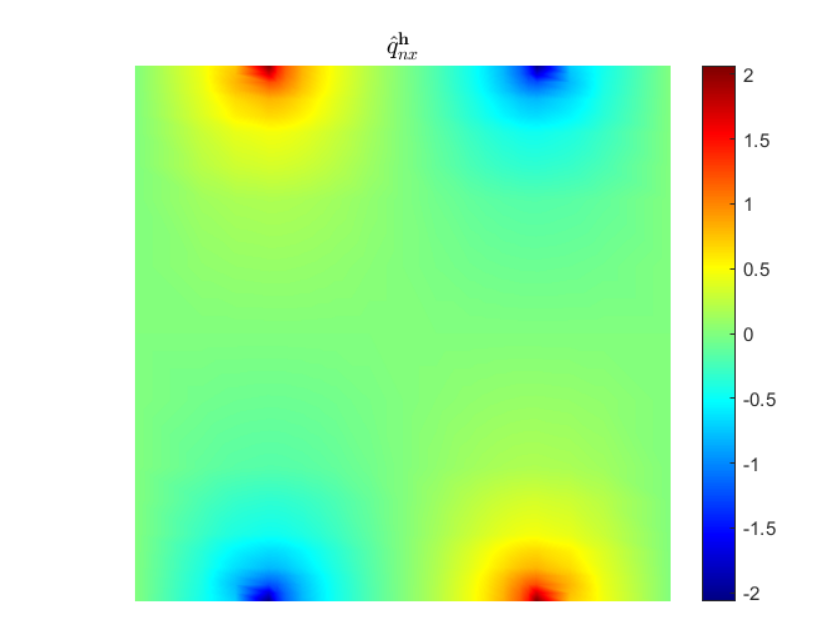}
    \end{subfigure}
    \begin{subfigure}{0.49\textwidth}
        \includegraphics[trim=3cm 10cm 4cm 9cm, clip=true, width=\textwidth]{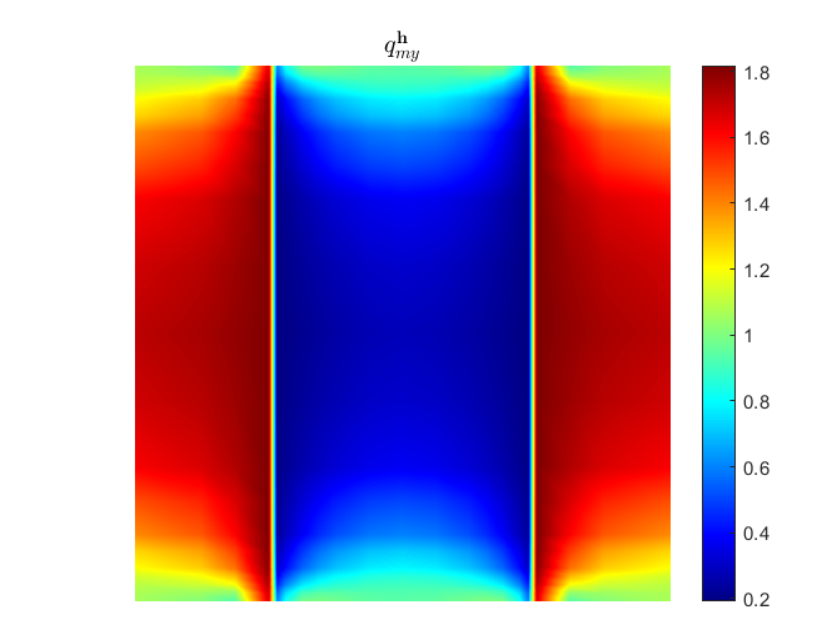}
    \end{subfigure}
    \begin{subfigure}{0.49\textwidth}
        \includegraphics[trim=3cm 10cm 4cm 9cm, clip=true, width=\textwidth]{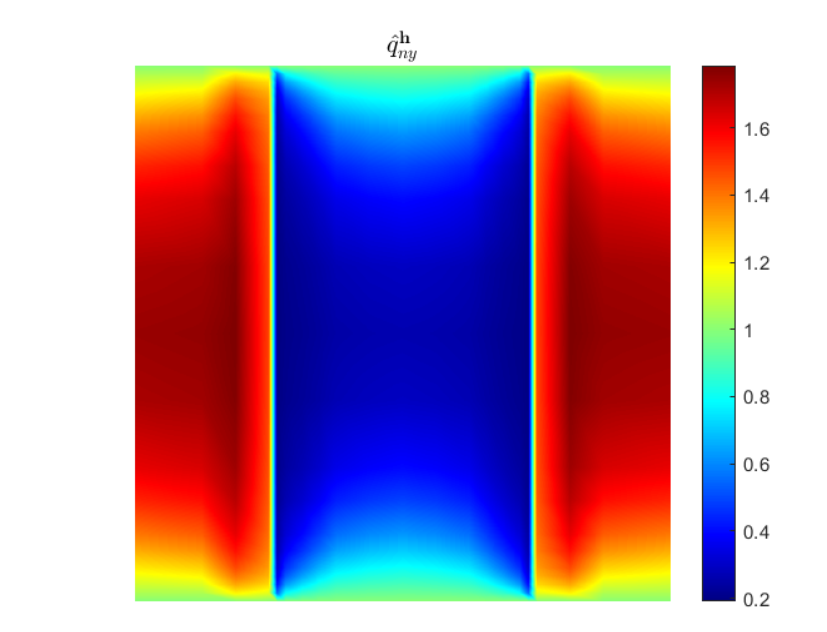}
    \end{subfigure}
    \caption{2D homogenization problem: components of the PGD flux $\bs{q}_m^{\bs{h}}$ (with 4 modes) and the SA PGD flux $\hat{\bs{q}}_n^{\bs{h}}$ (with 5 modes) at the end of the goal-oriented adaptive PGD procedure.}
    \label{fig:flux_pgd_hom}
\end{figure}

We use our goal-oriented adaptive PGD algorithm with two coarse initial meshes (in the $x$ and the $z$ direction) composed of 8 elements. The tolerance on $Q_2$ is set at 1\%. The PGD solution with 4 modes obtained at the end of the adaptive procedure is shown in Figure~\ref{fig:sol_pgd_hom}, along with the adapted meshes. The two components of the PGD flux and the SA PGD flux are shown in Figure~\ref{fig:flux_pgd_hom}.

\begin{figure}[h!]
    \centering
    \includegraphics[trim=3cm 9cm 4cm 9cm, clip=true, width=0.5\textwidth]{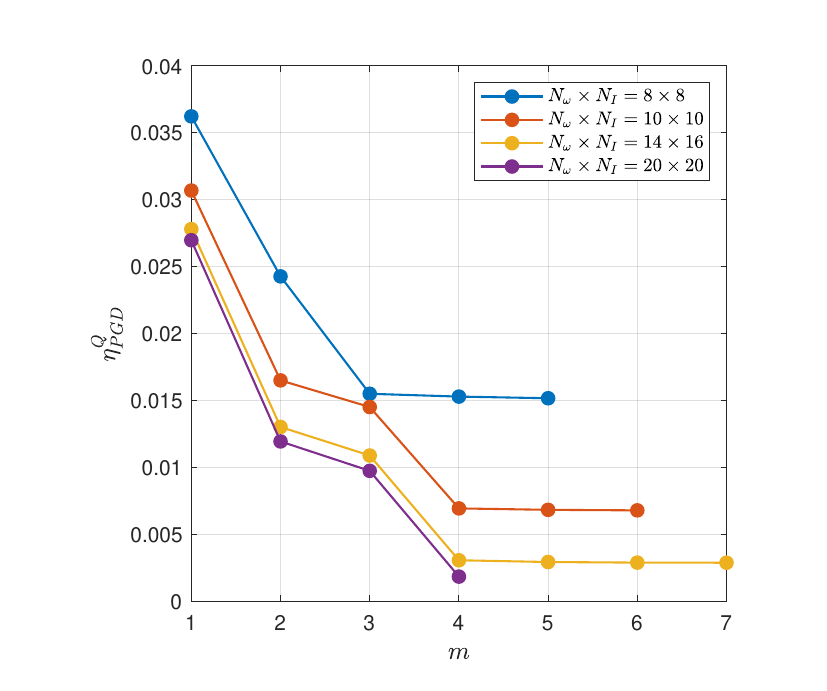}
    \caption{2D homogenization problem: convergence of the PGD error estimate for $Q_2$.}
    \label{fig:hist_conv_4.pdf}
\end{figure}

Eventually, the convergence of the error estimator on $Q_2$ is shown in Figure~\ref{fig:hist_conv_4.pdf}. We observe that for a given discretization, the error decreases monotonically before reaching a plateau. With regards to the accuracy of the estimator (with respect to the true error computed using a solution considered to be exact on a very fine mesh), we note that the effectivity indice is 1.74 at the end of the first mode adaptivity procedure and 1.38 at the end of the whole goal-oriented adaptive PGD procedure for $Q_2$.  

\section{Conclusion and perspectives} \label{sec:conclusion}

In this work, we addressed the certification of Proper Generalized Decomposition (PGD) reduced-order models based on the separation of spatial variables. Such representations are particularly attractive for the simulation of plate-like structures, since they allow the original problem to be reformulated as a sequence of lower-dimensional problems. A verification procedure based on the Constitutive Relation Error (CRE) concept was developed for this class of PGD approximations. The main difficulty lies in the construction of equilibrated fluxes compatible with the separated representation. To overcome this issue, a dedicated procedure was proposed, enabling the construction of statically admissible flux fields and, consequently, the derivation of guaranteed \textit{a posteriori} error bounds. The resulting estimator provides a rigorous upper bound on the total PGD error, accounting simultaneously for the discretization and reduction errors. Based on this estimator, an adaptive PGD strategy was introduced in order to automatically control both the discretization parameters and the number of PGD modes required to achieve a prescribed accuracy. The proposed verification framework was also extended to goal-oriented error estimation, allowing reliable control of quantities of interest. Numerical experiments conducted on two- and three-dimensional diffusion problems demonstrated the robustness of the proposed estimator, the quality of the associated error bounds, and the efficiency of the adaptive strategy. The present work demonstrates the possibility of certifying PGD approximations involving a full separation of spatial variables. Future developments will concern more complex classes of problems, such as linear elasticity, and curved geometries (shells).




\bibliographystyle{elsarticle-num} 
\bibliography{references}



\end{document}